\documentclass[11pt]{amsart}

\usepackage{a4wide}
\usepackage{bbm}
\usepackage{amsxtra} 
\usepackage{t1enc} 
\usepackage{color}
\usepackage{epsf}
\usepackage{amssymb}
\usepackage{graphicx}
\usepackage{stmaryrd}

\DeclareMathAlphabet{\mathpzc}{OT1}{pzc}{m}{it} %% little o

\newtheorem{theorem}{Theorem}[section]
\newtheorem{lem}[theorem]{Lemma}
\newtheorem{prop}[theorem]{Proposition}
\newtheorem{coro}[theorem]{Corollary}

\theoremstyle{definition}
\newtheorem{defi}[theorem]{Definition}
\newtheorem{rem}[theorem]{Remark}
\newtheorem{ex}[theorem]{Example}

\newcommand{\Z}{\mathbbm{Z}}
\newcommand{\N}{\mathbbm{N}}

\newcommand{\Q}{\mathbbm{Q}}

\newcommand{\R}{\mathbbm{R}}
\newcommand{\C}{\mathbbm{C}}

\newcommand{\Ha}{\mathbbm{H}}
\newcommand{\Ic}{\mathbbm{I}}

\begin{document}

\title[Discrete tomography of icosahedral model sets]{Discrete tomography
  of icosahedral model sets}

\author[C. Huck]{Christian Huck}
\address{\hspace*{-1em} Department of Mathematics and Statistics, 
  The Open University, Walton Hall, Milton Keynes, MK7 6AA, United Kingdom}
\email{c.huck@open.ac.uk}
\urladdr{http://www.mathematics.open.ac.uk/People/c.huck}

\begin{abstract} 
The discrete tomography of B-type and F-type icosahedral model
sets is investigated, with an emphasis on reconstruction and uniqueness problems. These are
motivated by the request of materials science for the unique reconstruction of
quasicrystalline structures from a small number of images produced by quantitative high
resolution transmission electron microscopy.
\end{abstract}

\maketitle

\section{Introduction}\label{intro}

{\em Discrete tomography} (the word ``tomography'' is derived from the Greek $\tau o\mu o\sigma$, meaning a slice) is concerned with the 
inverse problem of retrieving information about some {\em finite}
object from 
(generally noisy) information about its slices. A typical example is the {\em reconstruction} of a finite point set in Euclidean $3$-space 
from its line sums in a small number of directions. More precisely, a
({\em discrete parallel}\/) {\em
  X-ray} of a finite subset of
Euclidean $d$-space $\R^d$ in direction $u$ gives the
number of points of the set on each line in $\R^d$ parallel to
$u$. This concept should not be confused with $X$-rays in diffraction
theory, which provide rather different information on the
underlying structure that is based on statistical pair correlations;
compare~\cite{cow}, \cite{few} and \cite{gu}. In the
 classical setting, motivated by crystals, the positions to be determined form a subset of a common translate of the cubic lattice $\Z^3$ or, more
generally, of an arbitrary lattice $L$ in $\R^3$. In fact, many of the problems in discrete tomography
have been studied on $\Z^2$, the classical planar setting of
discrete tomography; see~\cite{HK}, \cite{Gr} and~\cite{GGP}. Beyond
the case of perfect crystals, one has to
take into account wider classes of sets, or at least
significant deviations from the lattice structure. As an intermediate
step between periodic and random (or amorphous) Delone sets, we consider systems of {\em aperiodic
  order}, more precisely, so-called {\em model sets} (or {\em
  mathematical quasicrystals}), which are commonly regarded as good mathematical models for quasicrystalline structures in
nature~\cite{St}. 

Our interest in the discrete tomography of
model sets is mainly motivated by the task of structure determination of
quasicrystals, a new type of solids discovered $25$ years ago;
see~\cite{s} for the pioneering paper and~\cite{so,j,d} for background and
applications. More precisely, we address the problem of uniquely reconstructing three-dimensional quasicrystals from their images under quantitative {\em high
  resolution transmission electron microscopy} (HRTEM) in a small
number of directions. In fact, in~\cite{ks} and ~\cite{sk} a technique is described, based
on HRTEM, which can effectively measure the number of atoms lying on
lines parallel to certain directions; it is called QUANTITEM
(\textbf{QU}antitative \textbf{AN}alysis of \textbf{T}he \textbf{I}nformation from \textbf{T}ransmission
\textbf{E}lectron \textbf{M}icroscopy). At present, the measurement of the number
of atoms lying on a line can only be approximately achieved for some
crystals; {\it cf.}~\cite{ks,sk}. However, it is reasonable to expect
that future developments in technology will improve this situation. 

In this text, we consider both {\em B-type} and {\em F-type icosahedral
model sets} $\varLambda$ in $3$-space which can be described in algebraic terms by using the {\em
icosian ring}; {\it cf.}~\cite{cmp}, ~\cite{Moody}
and~\cite{MP}. Note that the terminology originates from the fact that the
underlying $\Z$-modules (to be explained in Section~\ref{sec3}) of B-type and F-type icosahedral model sets
can be obtained as projections of body-centred and face-centred
hypercubic lattices in $6$-space, respectively. The F-type icosahedral phase is the most common among
the icosahedral quasicrystals. Below, we nevertheless develop the theory for both the B-type (also called I-type) and the F-type
phase. Well known examples of icosahedral
quasicrystals include the
aluminium alloys AlMn and AlCuFe; {\it cf.}~\cite{HG} for further examples.

In practice, only $X$-rays in
$\varLambda$-directions, {\it i.e.},
directions parallel to non-zero elements of the difference set
$\varLambda-\varLambda$ of $\varLambda$ ({\it i.e.}, the set of interpoint vectors of $\varLambda$) are reasonable. This is due to the fact that $X$-rays in non-$\varLambda$-directions are meaningless since the resolution coming from such $X$-rays would not be good enough to allow a quantitative analysis -- neighbouring lines are not sufficiently separated. In fact, in order to obtain applicable results, one even has to find $\varLambda$-directions that guarantee HRTEM images of high resolution, {\it i.e.}, yield dense lines in the corresponding quasicrystal $\varLambda$. 

Any lattice $L$ in $\R^d$ can be sliced into lattices of dimension $d-1$. More generally, model sets have a dimensional hierarchy, {\it i.e.}, any model set in
 $d$ dimensions can be sliced into model sets of dimension
 $d-1$. In Proposition~\ref{slice}, it is shown that {\it generic} (to
 be explained in Section~\ref{sec3}) B-type and F-type icosahedral model sets can be
 sliced into (planar) {\em cyclotomic model sets}, whose discrete tomography we have   
studied earlier; {\it cf.}~\cite{BG2,H2} and ~\cite{H}. The latter observation will be crucial,
since it enables us to use the results on the discrete tomography of
cyclotomic model sets, slice by slice.

Using the slicing of generic icosahedral model sets into cyclotomic
model sets and the results from~\cite{BG2}, it was shown in~\cite{H2} that the algorithmic 
problem of {\em reconstructing} finite subsets of a large class of
generic icosahedral model sets $\varLambda$ ({\it i.e.}, those with polyhedral windows) given $X$-rays in {\em two} $\varLambda$-directions can be solved in
polynomial time in the real RAM-model of computation (Theorem~\ref{th2}). Since this {\em reconstruction problem} can possess rather
different solutions, one is led to the investigation of the
corresponding {\em uniqueness problem}, {\it i.e.}, the (unique) {\em
  determination} of finite subsets of a fixed icosahedral model set
$\varLambda$ by $X$-rays in a
 small number of suitably prescribed $\varLambda$-directions. Here, a subset $\mathcal{E}$ of the set of all finite subsets of
a fixed icosahedral model set $\varLambda$ is said to be {\em determined} by the
$X$-rays in a finite set $U$ of directions if
different sets $F$ and $F'$ in $\mathcal{E}$ cannot have the same
$X$-rays in the directions of $U$. Since, as demonstrated in Proposition~\ref{source}, any fixed number of $X$-rays in
$\varLambda$-directions is insufficient to determine the entire class of
finite subsets of a fixed icosahedral model set $\varLambda$, it is necessary to impose some restriction in order to obtain positive uniqueness results. In Proposition~\ref{m+1}, it is shown that the finite subsets $F$ of cardinality less
than or equal to some $k\in \mathbbm{N}$ of a fixed icosahedral model
set $\varLambda$ are determined by any set of $k+1$ $X$-rays in
pairwise non-parallel $\varLambda$-directions. Proposition~\ref{bounded} then shows
 that, for every $R>0$ and any fixed
icosahedral model set $\varLambda$, there are two non-parallel $\varLambda$-directions such that the
 set of bounded subsets of $\varLambda$ with
diameter less than $R$ is determined by the $X$-rays in these
directions. For our main result, we restrict the set of finite subsets of a fixed
icosahedral model set $\varLambda$ by
considering the class of {\em convex} subsets of
$\varLambda$. They are finite sets $C\subset \varLambda$ whose convex
hulls contain no new points of $\varLambda$, {\it i.e.}, finite sets
$C\subset \varLambda$ with $C = \operatorname{conv}(C)\cap
\varLambda$. By using the slicing of generic icosahedral model sets into
cyclotomic model sets again, it is shown that there are {\em four} pairwise non-parallel $\varLambda$-directions such that the set of convex subsets of any
icosahedral model set $\varLambda$ are determined by their $X$-rays in these
directions (Theorem~\ref{main1ico}). In fact, it turns out that one can choose four $\varLambda$-directions which
provide uniqueness and yield dense lines in icosahedral model sets,
the latter making this result look promising in view of real
applications (Example~\ref{u5} and Remark~\ref{rempl}). Finally, we
demonstrate that, in an approximative sense, this result holds in a
far more general (and relevant) situation, where one deals with a whole family of
generic icosahedral model sets at the same time, rather than dealing with a single fixed icosahedral model set.   

\section{Preliminaries and notation}\label{found}

Natural numbers are always assumed to be positive, {\it i.e.}, 
$\mathbbm{N}=\{1,2,3,\dots\}$. Throughout the text, we use the
convention that the symbol $\,\subset\,$ includes equality. We denote the norm in Euclidean $d$-space $\mathbbm{R}^{d}$
by $\Arrowvert \cdot \Arrowvert$. The unit sphere in
$\mathbbm{R}^{d}$ is denoted by $\mathbb{S}^{d-1}$, {\it i.e.}, $\mathbb{S}^{d-1}=\{x\in\R^d\,|\,\Arrowvert x \Arrowvert =1\}$. Moreover, the elements of $\mathbb{S}^{d-1}$ are also called
{\em directions}. Recall that a {\em homothety} $h\!:\, \mathbbm{R}^{d} \rightarrow
\mathbbm{R}^{d}$ is given by $x \mapsto \lambda x + t$, where
$\lambda \in \R$ is positive and $t\in \mathbbm{R}^{d}$. We call a homothety {\em expansive} if $\lambda>1$. If $x\in\R$, then $\lfloor x \rfloor$
denotes the greatest integer less than or equal to $x$. For $r>0$ and $x\in\R^{d}$,
$B_{r}(x)$ is the open ball of radius $r$ about $x$. For a subset $S\subset
\mathbbm{R}^{d}$, $k\in \mathbbm{N}$ and $R>0$, we denote by
$\operatorname{card}(S)$, $\mathcal{F}(S)$, $\mathcal{F}_{\leq k}(S)$, $\mathcal{D}_{<R}(S)$,
$\operatorname{int}(S)$, $\operatorname{cl}(S)$, $\operatorname{bd}(S)$, $\operatorname{conv}(S)$, $\operatorname{diam}(S)$
and $\mathbbm{1}_{S}$ the cardinality, the set of finite subsets, the
set of finite subsets of $S$ having cardinality less than or equal to
$k$, the set of subsets of $S$
with diameter less than $R$, interior, closure,
boundary, convex hull, diameter and
characteristic function of $S$, respectively. The {\em
  centroid} (or {\em centre of mass}) of an element
$F\in\mathcal{F}(\R^d)$ is defined as
$(\sum_{f\in F}f)/\operatorname{card}(F)$. A
 linear subspace $T$ of $\mathbbm{R}^{d}$ is called an
$S${\em-subspace} if it is generated by elements of the {\em difference set} $S-S:=\{s-s'\,|\,s,s'\in S\}$ of $S$. A direction
$u\in\mathbb{S}^{d-1}$ is called an $S${\em-direction} if it is
parallel to a non-zero element of $S-S$. As usual, $R^{\times}$
denotes the group of units of a given ring $R$. Finally, for $(a,b,c)^t\in\R^3\setminus\{0\}$, we denote by $H^{(a,b,c)}$ the hyperplane in $\R^3$ orthogonal to $(a,b,c)^t$.

\begin{defi}\label{xray..}
Let $d\in \mathbbm{N}$ and let $F\in
\mathcal{F}(\mathbbm{R}^{d})$. Furthermore, let $u\in
\mathbb{S}^{d-1}$ be a direction and let $\mathcal{L}_{u}^{d}$ be the set
of lines in direction $u$ in $\mathbbm{R}^{d}$. Then, the ({\em
  discrete parallel}\/) {\em X-ray} of $F$ {\em in direction} $u$ is
the function $X_{u}F: \mathcal{L}_{u}^{d} \rightarrow
\mathbbm{N}_{0}:=\mathbbm{N} \cup\{0\}$, defined by $$X_{u}F(\ell) :=
\operatorname{card}(F \cap \ell\,) =\sum_{x\in \ell}
\mathbbm{1}_{F}(x)\,.$$ Moreover, the {\em support} $(X_{u}F)^{-1}(\N)$ of $X_{u}F$, {\it
  i.e.},
the set of lines in $\mathcal{L}_{u}^{d}$ which pass through at least one
point of $F$, is denoted by $\operatorname{supp}(X_{u}F)$. For $z\in
\mathbbm{R}^{d}$, we denote by $\ell_{u}^{z}$ the element of
$\mathcal{L}_{u}^{d}$ which passes through $z$. Moreover, for $S\subset \mathbbm{R}^{d}$, we denote
by $\mathcal{L}_{u}^{S}$ the subset of $\mathcal{L}_{u}^{d}$ consisting of
all elements of the form $\ell_{u}^{z}$, where $z\in S$, {\it
  i.e.},
 lines in $\mathcal{L}_{u}^{d}$ which pass through at
least one point of $S$.
\end{defi}

\begin{lem}{\rm \cite[Lemma 5.1 and Lemma~5.4]{GG}}\label{cardinality}
Let $d\in \mathbbm{N}$ and let $u\in \mathbb{S}^{d-1}$ be a
direction. For all $F, F'\in \mathcal{F}(\mathbbm{R}^{d})$, one has: 
\begin{itemize}
\item[{\rm (a)}]
$X_{u}F=X_{u}F'$ implies $\operatorname{card}(F)=\operatorname{card}(F')$.
\item[{\rm (b)}]
If $X_{u}F=X_{u}F'$, the centroids of $F$ and $F'$ lie on the same line parallel to $u$.
\end{itemize}
\end{lem}

\begin{defi}\label{defgrid}
Let $d\geq 2$, let $U\subset\mathbb{S}^{d-1}$
be a finite set of pairwise
non-parallel directions and let $F\in\mathcal{F}(\R^d)$. We define the \emph{grid} of $F$ with respect to the $X$-rays in the
directions of $U$ as
$$
G^{F}_{U}\,\,:=\,\,\bigcap_{u\in U}\,\,\left( \bigcup_{\ell \in
  \mathrm{supp}(X_{u}F)} \ell\right)\,.
$$
\end{defi}

The following property follows immediately from the definition of grids.

\begin{lem}\label{fgrid}
Let $d\geq 2$. If $U\subset\mathbb{S}^{d-1}$
is a finite set of pairwise
non-parallel directions, then for all $F,F'\in\mathcal{F}(\R^d)$, one has
$$
(X_{u}F=X_{u}F'\;\,\forall u \in U) \;  \Longrightarrow\; F,F'\,\,\subset\,\,G^{F}_{U}\,\,=\,\,G^{F'}_{U}\,.
$$

\end{lem}

\begin{defi}
Let $d\geq 2$, let $\mathcal{E}\subset
\mathcal{F}(\mathbbm{R}^{d})$, and let $m\in\N$. Further, let $U\subset\mathbb{S}^{d-1}$ be a
finite set of directions. We say that $\mathcal{E}$ is {\em determined} by the $X$-rays in the directions of $U$ if, for all $F,F' \in \mathcal{E}$, one has
$$
(X_{u}F=X_{u}F'\;\,\forall u \in U) \;  \Longrightarrow\; F=F'\,.
$$
Further, we say that $\mathcal{E}$ is {\em determined} by $m$ $X$-rays if there exists a set $U$ of $m$ pairwise non-parallel directions such that $\mathcal{E}$ is determined by the $X$-rays in the directions of $U$. 
\end{defi}

The following property is straight-forward.

\begin{lem}\label{homotu}
Let $d\geq 2$, let $h\!:\,\mathbbm{R}^{d}  \rightarrow \mathbbm{R}^{d}$ be a
  homothety, and let $U\subset \mathbb{S}^{d-1}$ be a finite set of
  directions. Then, if $F$ and $F'$ are elements of $\mathcal{F}(\mathbbm{R}^d)$
with the same $X$-rays in the directions of $U$, the images $h(F)$ and $h(F')$ also have the same $X$-rays in the directions of $U$.
\end{lem}

Gardner and Gritzmann introduced the so-called \emph{convex lattice sets}, {\it i.e.}, finite subsets $C$ of some lattice $L\subset\R^d$ with $C =
\operatorname{conv}(C)\cap L$; {\it cf.}~\cite[Section 2]{GG}. More generally, we define as follows.

\begin{defi}\label{deficonvex} 
Let $d\in \N$ and let $S \subset \R^d$.
 A finite subset $C$ of $S$ is called a {\em convex subset of} $S$ if  it satisfies the equation  
$C =
\operatorname{conv}(C)\cap S$. Moreover, the set of all convex subsets of
$S$ is denoted by $\mathcal{C}(S)$. 
\end{defi}

\section{Icosahedral model sets}\label{sec3}

We shall always denote the
golden ratio by $\tau$, {\it i.e.}, $\tau=(1+\sqrt{5})/2$. Moreover, by $.'$ we will denote the unique
non-trivial Galois automorphism of the real quadratic number field $\Q(\tau)=\Q(\sqrt{5})=\Q\oplus\Q\tau$ (determined by $\sqrt{5}\mapsto
-\sqrt{5}$), whence $\tau'=-1/\tau=1-\tau$. Note that $\tau$ is an
algebraic integer (a root of $X^2-X-1\in \Z[X]$) of degree $2$ over
$\Q$. Moreover, $
\Z[\tau]=\Z\oplus\Z\tau$ is the ring of integers in $\Q(\tau)$ and, for its group of units, one further has $\Z[\tau]^{\times}=\{\tau^s\,|\,s\in\Z\}$ ({\it i.e.}, $\tau$ is a fundamental unit of $\Z[\tau]$); {\it cf.}~\cite{hw}.

\subsection{Definition and properties of icosahedral model sets}\label{sec31}

Let $\Ha$ be the skew field of {\em Hamiltonian quaternions}, {\it i.e.},
$$
\Ha =\{a+bi+cj+dk\,|\,a,b,c,d\in\R\}\,,$$
a four-dimensional vector space over $\R$ with a non-commutative
multiplication determined by the following relations for the generating elements $1$ (implicit in the above representation) and $i,j,k$:
$$
i^2=j^2=k^2=ijk=-1\,,
$$
together with the requirement that $\R$ is central in $\Ha$. Note that $\R$ is precisely the center of $\Ha$. The {\em conjugate} of $\alpha=a+bi+cj+dk\in\Ha$ is
defined by $\bar{\alpha}=a-bi-cj-dk$, the {\em reduced norm} by
$\operatorname{nr}(\alpha)=\alpha \bar{\alpha}=a^2+b^2+c^2+d^2$ and the {\em reduced trace}
 by $\operatorname{tr}(\alpha)=\alpha+\bar{\alpha}=2a$. Moreover, we shall
 sometimes call $\operatorname{Re}(\alpha):=a\in\R$ the {\em real part} and
 $\operatorname{Im}(\alpha):=(b,c,d)^t\in\R^3$ the {\em imaginary part} of
 $\alpha$. Let $
\Ha_{0}$ be the set of quaternions with real part $0$, {\it i.e.}, 
$$
\Ha_{0}
:=\{\alpha\in\Ha\,|\,\operatorname{tr}(\alpha)=0\}=\{bi+cj+dk\,|\,b,c,d\in\R\}\,\,\simeq\,\,
\R^3\,.
$$
The {\em icosian ring} $\Ic$ ({\em cf.}~\cite{cmp,Moody,MP}) is the additive
subgroup of $\Ha$ that is given by the integer linear combinations of the quaternions
$$
\left((\pm1,0,0,0)^{t}\right)^{\mathsf{A}},\tfrac{1}{2}\left((\pm1,\pm1,\pm1,\pm1)^{t}\right)^{\mathsf{A}},
\tfrac{1}{2}\left((0,\pm1,\pm\tau',\tau)^{t}\right)^{\mathsf{A}}\,,
$$
where we identify $\Ha$ with $\R^4$ via the basis $\{1,i,j,k\}$ and,
as in~\cite[Chapter 8]{con},
 the superscript $\sf{A}$ indicates that all even permutations of the
coordinates are allowed. The members of $\Ic$ are called {\em icosians}. Note that $\Ic$ is a
ring, because these generators (which have reduced norm $1$) form a
multiplicative group, the {\em icosian group}, of order $120$. Note further that $\Ic$ is also a free $\Z[\tau]$-module of
rank $4$. By~\cite{BPR}, $\Ic$ is a maximal order of the quaternion algebra
$\Ha(\mathbbm{Q}(\tau))$ over $\mathbbm{Q}(\tau)$, defined similar to $\Ha$ as
$$
\Ha\big(\mathbbm{Q}(\tau)\big) =\big\{a+bi+cj+dk\,\big|\,a,b,c,d\in\mathbbm{Q}(\tau)\big\}\,.
$$
The set
$$
\Ic_{0}\,:=\,\operatorname{Im}(\Ic \,\cap \,\Ha_{0})\,\subset\,\R^3
$$
of `pure imaginary' icosians is generated as an additive group by the elements
$$
\left((\pm1,0,0)^{t}\right)^{\mathsf{A}},\tfrac{1}{2}\left((\pm1,\pm \tau',\pm \tau)^{t}\right)^{\mathsf{A}}\,,$$
where the superscript $\sf{A}$ is defined as above. Consider the standard
body-centred icosahedral module $\mathcal{M}_{\text{B}}$ of
quasicrystallography, defined as 
\begin{equation}\label{mbeqn}
\begin{split}
&\mathcal{M}_{\text{B}}\,\,:=\,\,\Z[\tau](2,0,0)^{t}\oplus\Z[\tau](1,1,1)^{t}\oplus\Z[\tau](\tau,0,1)^{t}\\[2mm]
&\,\,\,\,\,\,\,\,\,\,\,\,\,\,=\,\,\Z[\tau](0,2,0)^{t}\oplus\Z[\tau](-1,-\tau',\tau)^{t}\oplus\Z[\tau](1,1,1)^{t}\\[2mm]
&\,\,\,\,\,\,\,\,\,\,\,\,\,\,=\,\,\left\{(\beta,\gamma,\delta)^{t}\,\left
  |\,\begin{array}{l}\beta,\gamma,\delta\in\Z[\tau]\mbox{, with}\\
    \tau^2\beta+\tau\gamma+\delta\equiv\; 0 \;(\operatorname{mod} 2)
\end{array}\right\}\,;\right.
\end{split}
\end{equation}
{\it cf.}~\cite{baake,BPR} and references therein. One has
$
\operatorname{Im}(\Ic)=\tfrac{1}{2}\mathcal{M}_{\text{B}}
$ and, further, 
$
\Ic_0=\tfrac{1}{2}\mathcal{M}_{\text{F}}
$, where $\mathcal{M}_{\text{F}}$ is the standard face-centred icosahedral module of
quasicrystallography, defined as 
\begin{equation}\label{mbeqn2}
\begin{split}
&\mathcal{M}_{\text{F}}\,\,:=\,\,\left\{(\beta,\gamma,\delta)^{t}\,\left
  |\,\begin{array}{l}\beta,\gamma,\delta\in\Z[\tau]\mbox{, with}\\
    \beta\equiv\tau\gamma\equiv\tau^2\delta\;(\operatorname{mod} 2)
\end{array}\right\}\right.
\\[2mm]
&\,\,\,\,\,\,\,\,\,\,\,\,\,\,=\,\,\left\{\left.(\beta,\gamma,\delta)^{t}\in\mathcal{M}_{\text{B}}\,\right
  |\,\beta+\gamma+\delta\equiv\; 0 \;(\operatorname{mod} 2)\right\}
\\[2mm]
&\,\,\,\,\,\,\,\,\,\,\,\,\,\,=\,\,\Z[\tau](2,0,0)^{t}\oplus\Z[\tau](\tau+1,\tau,1)^{t}\oplus\Z[\tau](0,0,2)^{t}
\\[2mm]
&\,\,\,\,\,\,\,\,\,\,\,\,\,\,=\,\,\Z[\tau](0,2,0)^{t}\oplus\Z[\tau](-1,-\tau',\tau)^{t}\oplus\Z[\tau](2,0,0)^{t}\,\,\stackrel{4}{\subset}\,\,\mathcal{M}_{\text{B}}\,,
\end{split}
\end{equation}
where integers on top of the inclusion symbol denote the corresponding subgroup indices; {\it cf.}~\cite{baake,BPR} again. Both $\mathcal{M}_{\text{B}}$ and $\mathcal{M}_{\text{F}}$ are free
$\Z[\tau]$-modules of rank $3$, and are hence $\Z$-modules of rank
$6$. Moreover, both $\mathcal{M}_{\text{B}}$ and $\mathcal{M}_{\text{F}}$ have icosahedral symmetry, {\it i.e.},
they are 
invariant under the action of the rotation group $Y$. This group is generated by the
rotations which are given, with respect to the canonical basis, by the following matrices
\begin{eqnarray}\label{rep}
\left( \begin{array}{ccc}-1&0&0\\0&-1&0\\0&0&1\end{array}
\right)
\,,\,\,\,\, \frac{1}{2}\left( \begin{array}{ccc}\tau&-1&-\tau'\\1&-\tau'&-\tau\\-\tau'&\tau&1\end{array}
\right)
\,.
\end{eqnarray}
Note that $Y$ is the rotation group of the regular icosahedron centred at the origin $0\in\R^3$ with orientation such that each
coordinate axis passes through the mid-point of an edge, thus
coinciding with $2$-fold axes of the icosahedron. Moreover, the
matrix on the left (resp., on the right) is an order $2$ (resp., order
5)
rotation. 

\begin{rem}
There is
another $\Z$-module of rank $6$, intermediate between
$\mathcal{M}_{\text{F}}$ and $\mathcal{M}_{\text{B}}$, which also has
icosahedral symmetry. This is the standard primitive icosahedral module $\mathcal{M}_{\text{P}}$, defined as
$$
\mathcal{M}_{\text{P}}\,\,:=\,\,\left\{\left.(\beta,\gamma,\delta)^{t}\in\mathcal{M}_{\text{B}}\,\right
  |\,\beta+\gamma+\delta\equiv\; 0 \mbox{ or } \tau
  \;(\operatorname{mod} 2)\right\}\,.$$
In contrast to $\mathcal{M}_{\text{F}}$ and $\mathcal{M}_{\text{B}}$,
$\mathcal{M}_{\text{P}}$ fails to be a $\Z[\tau]$-module. In fact,
$\mathcal{M}_{\text{P}}$ is a $\Z[2\tau]$-module only, and it is a $\Z$-module of rank $6$.  
\end{rem}

By definition, {\em model sets} arise from so-called {\em cut and
  project schemes}; {\it cf.}~\cite{BM,Moody} for general background
material and see~\cite{B} for a gentle introduction. In the case of
 Euclidean internal spaces, these are
  commutative diagrams of the following form, where $\pi$ and $\pi_{\textnormal{\tiny int}}$ denote
the canonical projections; {\it cf.}~\cite{Moody}.  

\begin{equation}\label{cutproj}
\renewcommand{\arraystretch}{1.2}
\begin{array}{ccccc}
\,\,\,& \pi & & \pi_{\textnormal{\tiny int}} & \vspace*{-1.5ex} \\
\mathbbm{R}^{d}\,\,\, & \longleftarrow & \mathbbm{R}^{d}\times \mathbbm{R}^{m}  & \longrightarrow & \mathbbm{R}^{m} \\
\cup\,\,\,\,\,&&\,\,\,\,\,\,\,\,\,\,\,\,\cup\mbox{\tiny\, lattice}&&\,\,\,\,\,\,\,\,\cup\mbox{\tiny\, dense}\\
 \,\,\,& \mbox{\tiny 1--1} & &  & \vspace*{-1.5ex} \\
 \!\!\!L \,\,\,& \longleftrightarrow & \!\!\widetilde{L}& \longrightarrow &\!\!L^{\star} \\
\end{array}
\end{equation}

Here, $\tilde{L}$ is a lattice in $\R^d\times \R^m$. Further, we
assume that 
the restriction $\pi|_{\tilde{L}}$ is injective and that the image $\pi_{\textnormal{\tiny
    int}}(\tilde{L})$ is a dense subset of $\R^m$. Letting
$L:=\pi(\tilde{L})$, the bijectivity of the (co-)restriction
$\pi|_{\tilde{L}}^L$ allows us to define a map $.^{\star}\! : \,
L\rightarrow \R^m$ by $\alpha^{\star}:=\pi_{\textnormal{\tiny
    int}}({(\pi|_{\tilde{L}}^L)}^{-1}(\alpha))$. Then, one has $L^{\star}=\pi_{\textnormal{\tiny
    int}}(\tilde{L})$ and, further,
$\tilde{L}=\{(l,l^{\star})\,|\,l\in L\}$.

\begin{defi}\label{mset}
Given a subset $W\subset \R^m$ with
$\varnothing\, \neq\, \operatorname{int}(W)\subset W\subset \operatorname{cl}(\operatorname{int}(W))$
and $\operatorname{cl}(\operatorname{int}(W))$ compact, a so-called {\em window}, and any $t\in\R^d$, we obtain a model set $$\varLambda(t,W) := t+\varLambda(W)$$ relative
to the above cut and project scheme~(\ref{cutproj}) by setting
$$\varLambda(W):=\{\alpha\in L\,|\,\alpha^{\star}\in W\}\,.$$
Moreover, $\mathbbm{R}^{d}$ (resp., $\mathbbm{R}^{m}$) is called the
\emph{physical} (resp., \emph{internal}) space. The map
$.^{\star}\! : \,  L\rightarrow \R^m$, as defined above, is the
so-called \emph{star map} of $\varLambda(t,W)$, $W$ is referred to
as the {\em window} of $\varLambda(t,W)$ and $L$ is the so-called 
\emph{underlying $\Z$-module of $\varLambda(t,W)$}. The model set
$\varLambda(t,W)$ is called {\em generic} if it satisfies $\operatorname{bd}(W)\cap
L^{\star}=\varnothing$. Moreover, it is called {\em regular} if the boundary $\operatorname{bd}(W)$ has Lebesgue
measure $0$ in $\R^m$.
\end{defi}

\begin{rem}
Every translate of a window $W\subset
\R^m$ is a window again.
\end{rem}

\begin{defi}\label{defibf}
{\em B-type icosahedral model sets} $\varLambda^{\rm
  B}_{{\rm ico}}(t,W)$ arise from the cut and
project scheme~\eqref{cutproj} by setting $d:=m:=3$, 
$L:=\operatorname{Im}(\Ic)$ and letting the star map $
.^{\star}\!:\, \operatorname{Im}(\Ic)
\rightarrow \R^3
$
be defined by applying the Galois conjugation $.'$ to each
coordinate of an element $\alpha\in \operatorname{Im}(\Ic)$. We denote by $\mathcal{I}^{\rm B}$ the set of all B-type icosahedral model sets
and define $\mathcal{I}^{\rm B}_g$ as the subset  of
all generic B-type icosahedral model sets. Additionally, for a window $W\subset
\R^3$, we set
$$
\mathcal{I}^{\rm B}_{g}(W)\,\,:=\,\,\{\varLambda^{\rm B}_{{\rm
    ico}}(t,s + W)\,|\,t,s\in\mathbbm{R}^{3} \}\,\,\cap\,\,
\mathcal{I}^{\rm B}_g \,.
$$ 
{\em F-type icosahedral model sets} $\varLambda^{\rm
  F}_{{\rm ico}}(t,W)$ arise from the cut and
project scheme~\eqref{cutproj} by setting $d:=m:=3$, 
$L:=\Ic_0$ and letting the star map $
.^{\star}\!:\, \Ic_0
\rightarrow \R^3
$
again be defined by applying the Galois conjugation $.'$ to each
coordinate of an element $\alpha\in \Ic_0$. Moreover, the
sets $\mathcal{I}^{\rm F}$, $\mathcal{I}^{\rm F}_g$ and $
\mathcal{I}^{\rm F}_{g}(W)$, where $W\subset
\R^3$ is a window, are defined similarly. Below, we say that $\varLambda_{{\rm ico}}(t,W)$ 
is an {\em icosahedral model set} if $\varLambda_{{\rm
    ico}}(t,W)=\varLambda_{{\rm ico}}^{\rm B}(t,W)$ or $\varLambda_{{\rm
    ico}}(t,W)=\varLambda_{{\rm ico}}^{\rm F}(t,W)$. Finally, B-type (resp., F-type) icosahedral model
sets are also referred to as \emph{icosahedral model sets with underlying 
$\Z$-module $\operatorname{Im}(\Ic)$} (resp., $\Ic_0$).  
\end{defi}

\begin{rem}
Both star maps as defined in
Definition~\ref{defibf} are $\Q$-linear monomorphism of Abelian groups
and naturally extend to a monomorphism $\Q(\tau)^3
\rightarrow \R^3
$, which we also denote by $.^{\star}$. Both in the B-type and the
F-type case, we shall denote by $.^{-\star}$
    the inverse of the co-restriction of the corresponding star map $.^{\star}\!:\, L
\rightarrow L^{\star}$ to
its image. The images of both maps $\,\widetilde{.}\!:\, L
\rightarrow \R^3\times \R^3$, defined by $\alpha\mapsto
(\alpha,\alpha^{\star})$, are indeed lattices in
$\R^3\times \R^3\simeq\R^6$. In fact, these images have a natural interpretation as a weight lattice
of type $D_6^{\ast}$ in the B-type case and a root lattice
of type $D_6$ in the F-type case; {\it cf.}~\cite{cmp,con} for
background. Finally, one can easily verify that the images $\operatorname{Im}(\Ic)^{\star}$ and
${\Ic_0}^{\star}$ are indeed dense subsets of $\R^3$.
\end{rem}

We refer the reader to~\cite{Moody,PABP} for details and related general settings, and to~\cite{BM} for general background. Before we collect some properties of
icosahedral model sets, recall the following notions. A subset $\varLambda$ of $\R^d$ is called {\em uniformly discrete} if there is a radius
$r>0$ such that every ball $B_{r}(x)$ with $x\in\mathbbm{R}^{d}$ contains at most one point of
  $\varLambda$. Further, $\varLambda$ is called {\em relatively dense} if there is a radius $R>0$
  such that every ball $B_{R}(x)$ with 
  $x\in\mathbbm{R}^{d}$ contains at least one point of $\varLambda$.

\begin{rem}\label{propcm}
Let $\varLambda$ be an icosahedral model set with window $W$. Then, $\varLambda$ is a {\em Delone set} in $\R^3$ ({\it i.e.}, $\varLambda$ is both uniformly discrete and relatively dense) and
is of {\em finite local complexity} ({\it i.e.},
 $\varLambda-\varLambda$ is closed and discrete). Note that $\varLambda$ is of finite local complexity if and
  only if for every $r>0$ there
  are, up to translation, only finitely many point sets (called
  \emph{patches of diameter} $r$) of the form $\varLambda\cap
  B_{r}(x)$, where $x\in\mathbbm{R}^3$; {\it cf.}~\cite[Proposition 2.3]{Schl}. In fact, $\varLambda$ is even a {\em
  Meyer set}, {\it i.e.}, $\varLambda$ is a Delone set and $\varLambda-\varLambda$ is uniformly
discrete; compare \cite{Moody}. Further, $\varLambda$ is an {\em
  aperiodic} model set, {\it i.e.}, $\varLambda$ has no translational
symmetries. Moreover, if
$\varLambda$ is {\em regular}, $\varLambda$ is {\em
  pure point diffractive}, {\it i.e.}, the Fourier transform of the autocorrelation density
that arises by placing a delta peak (point mass) on each point of $\varLambda$ looks purely point-like; {\it cf.}~\cite{Schl}. If $\varLambda$ is
  generic, $\varLambda$ is {\em repetitive}, {\it i.e.}, given any patch of radius $r$, there is a radius $R>0$ such that any ball of radius $R$ contains at least one translate of this patch; {\it cf.}~\cite{Schl}. If $\varLambda$ is regular, the frequency of repetition of finite
patches is well defined, {\it i.e.}, for any patch of radius $r$, the number of occurrences of translates of this patch per unit volume in the ball $B_r(0)$ of radius $r>0$ about the origin $0$ approaches a non-negative limit as $r\rightarrow \infty$; {\it cf.}~\cite{Schl2}. Moreover, if
$\varLambda$ is both generic and regular, and, if a suitable translate
of the window $W$ 
 has full icosahedral symmetry ({\it i.e.}, if a suitable translate
of the window $W$ is invariant under the action of the group
$Y_h^{\star}$ of order $120$, where
$Y_h^{\star}:=Y^{\star}\,\cup (-Y^{\star})$ and $Y^{\star}$ is the group
of rotations of order $60$ generated by
 the two matrices that arise from the two matrices in~\eqref{rep} by applying the conjugation $.'$ to each entry), then $\varLambda$ 
 has full icosahedral symmetry $Y_h:=Y\, \cup (-Y)$ in the sense of symmetries of LI-classes,
meaning that a discrete structure has a certain symmetry if the
original and the transformed structure are locally indistinguishable
(LI) ({\it i.e.}, up to translation, every {\em finite} patch in $\varLambda$ also appears in 
any of the other elements of its LI-class and {\em vice versa}); see~\cite{B} for details. Typical examples are balls and suitably oriented
versions of the icosahedron, the dodecahedron, the rhombic
triacontahedron (the latter also known as Kepler's body) and its dual,
the icosidodecahedron. 
\end{rem}

\begin{figure}
\centerline{\epsfysize=0.62\textwidth\epsfbox{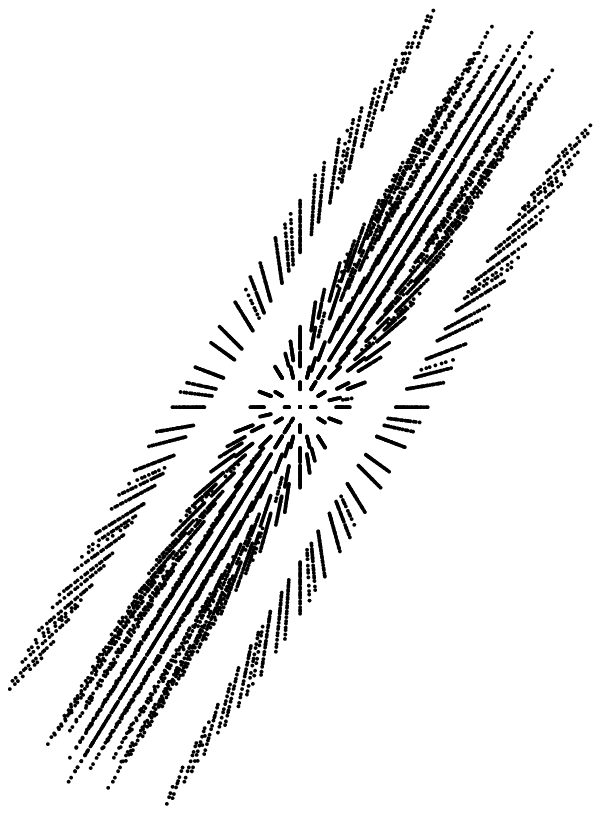}\hspace{0.00\textwidth}
\epsfysize=0.62\textwidth\epsfbox{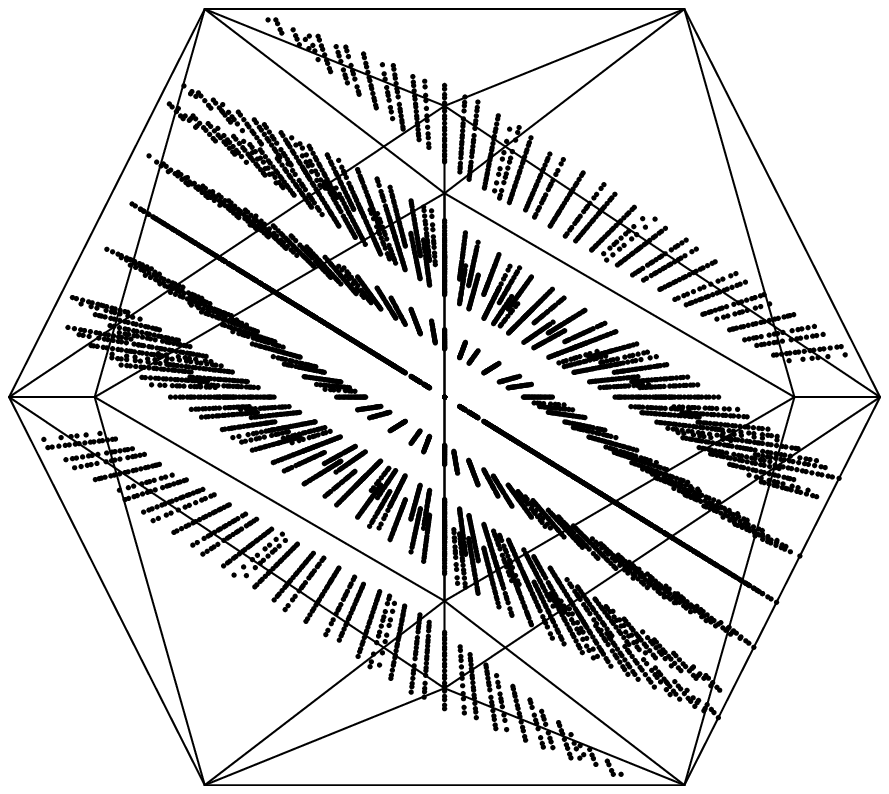}}
\caption{A few slices of a patch of the icosahedral model set
  $\varLambda_{\rm ico}^{\rm B}$ (left) and their $.^{\star}$-images inside the icosahedral window in the internal space (right), both seen from the positive $x$-axis.}
\label{fig:slices}
\end{figure}

\begin{ex}\label{lambdaico}
For a generic regular icosahedral model set with full icosahedral symmetry $Y_h$, consider $\varLambda_{\rm ico}^{\rm B}:=\varLambda_{\rm ico}^{\rm B}(0,s+W)$, where $s:=10^{-3}(1,1,1)^t$ and $W$ is the regular icosahedron with vertex set $Y_h^{\star}(\tau',0,1)^t$; see Figure~\ref{fig:slices} for an illustration.     
\end{ex}

\subsection{Cyclotomic model sets as planar sections of icosahedral model sets}\label{cycsec}

In this section, we shall demonstrate that both B-type and F-type icosahedral model sets $\varLambda$ can be
nicely sliced 
into {\em cyclotomic model sets with underlying $\Z$-module $\Z[\zeta_5]$}, where the slices
are intersections of $\varLambda$ with translates of the hyperplane $H^{(\tau,0,1)}$ in $\R^3$ orthogonal
to $(\tau,0,1)^t$. From now on, we always let $\zeta_{5}:=e^{2\pi i/5}$, as a specific choice of a primitive $5$th root of
unity in $\C$. Occasionally, we 
identify $\C$ with $\R^{2}$.

\begin{rem}\label{remzm2}
It is well known that the $5$th cyclotomic field $\Q(\zeta_{5})$ is an
algebraic number field of degree $4$ over $\Q$. Moreover, the field extension $\Q(\zeta_5)/ \mathbbm{Q}$
is a Galois extension with Abelian Galois group $G(\Q(\zeta_5)/
\mathbbm{Q}) \simeq (\Z / 5\Z)^{\times}$,
where $a\, (\textnormal{mod}\, 5)$ corresponds to the automorphism
given by\/ $\zeta_5 \mapsto \zeta_5^{a}$; {\it cf.}~\cite[Theorem
  2.5]{Wa}. Note that, restricted to the quadratic field $\Q(\tau)$,
both the Galois automorphism of $\Q(\zeta_5)/ \mathbbm{Q}$ that is
given by $\zeta_5 \mapsto \zeta_5^{3}$ and its complex conjugate
automorphism ({\it i.e.}, the automorphism given by $\zeta_5 \mapsto
\zeta_5^{2}$) induce the unique non-trivial Galois automorphism $.'$ of $\Q(\tau)/\Q$ (determined by $\tau\mapsto 1-\tau$). Further, $\Z[\zeta_{5}]$ is the
ring of integers in $\Q(\zeta_{5})$; {\it cf.}~\cite[Theorem 2.6]{Wa}. The ring
$\Z[\zeta_{5}]$ also is a $\Z[\tau]$-module of rank two. More
precisely, one has the equality $\Z[\zeta_{5}]=\Z[\tau]\oplus\Z[\tau]\zeta_{5}$; {\it cf.}~\cite[Lemma 1(a)]{BG2}. Since $\zeta_{5}^{3}$
is also a
primitive $5$th root of unity in $\C$, one further has the equality $\Z[\zeta_{5}]=\Z[\zeta_{5}^{3}]=\Z[\tau]\oplus\Z[\tau]\zeta_{5}^{3}$. 
\end{rem}

\begin{defi}\label{deficyc}
{\em Cyclotomic model sets with underlying $\Z$-module $\Z[\zeta_5]$}
$\varLambda_{{\rm cyc}}(t,W)$ arise from the cut and
project scheme~\eqref{cutproj} by setting $d:=m:=2$, $L:=\Z[\zeta_5]$ and letting
the star map $.^{\star_{5}}\!:\,L\rightarrow \R^2$ be either given by the non-trivial Galois automorphism of $\mathbbm{Q}(\zeta_5)/\Q$, defined by $\zeta_{5}\mapsto
\zeta_{5}^{3}$, or its complex conjugate automorphism. 
\end{defi}

\begin{rem}
The star map $.^{\star_{5}}$ as defined in Definition~\ref{deficyc} is
a monomorphism of Abelian groups. Further, the
image of the map $\,\tilde{.}^{_{5}}\!:\,L\rightarrow \R^2\times\R^2$,
defined by $\alpha\mapsto (\alpha,\alpha^{\star_{5}})$, is indeed a lattice in $\R^2\times\R^2$. Finally,
one can verify that the image $L^{\star_{5}}$ is indeed a dense subset of $\R^2$.
%\begin{defi}\label{deficyc}
%Given any subset $W\subset \R^2$ with
%$\varnothing\, \neq\, \operatorname{int}(W)\subset W\subset \operatorname{cl}(\operatorname{int}(W))$
%and $\operatorname{cl}(\operatorname{int}(W))$ compact, a so-called {\em window}, and any $t\in\R^2$, we obtain a
%planar model set, a so-called {\em
%  cyclotomic model set with underlying $\Z$-module $\Z[\zeta_5]$}, $$\varLambda_{{\rm cyc}}(t,W) := t+\varLambda_{{\rm cyc}}(W)$$ by setting
%$$\varLambda_{{\rm cyc}}(W):=\{z\in\Z[\zeta_5]\,|\,z^{\star_5}\in W\}\,,$$
%where the {\em star map} $.^{\star_{5}}$ of $\varLambda_{{\rm cyc}}(t,W)$ is either the non-trivial Galois automorphism of $\mathbbm{Q}(\zeta_5)/\Q$, defined by $\zeta_{5}\mapsto
%\zeta_{5}^{3}$, or its complex conjugate automorphism.
%\end{defi}
For the general setting, we refer the reader
to~\cite{BG2,H2,H}. By~\cite[Lemma 1.84(a)]{H2} (see also~\cite[Lemma
  25(a)]{H}), for all cyclotomic model sets $\varLambda$ with underlying
$\Z$-module $\Z[\zeta_5]$, the set of $\varLambda$-directions is precisely the set of $\Z[\zeta_5]$-directions.  
\end{rem}

\begin{figure}
\centerline{\epsfysize=0.48\textwidth\epsfbox{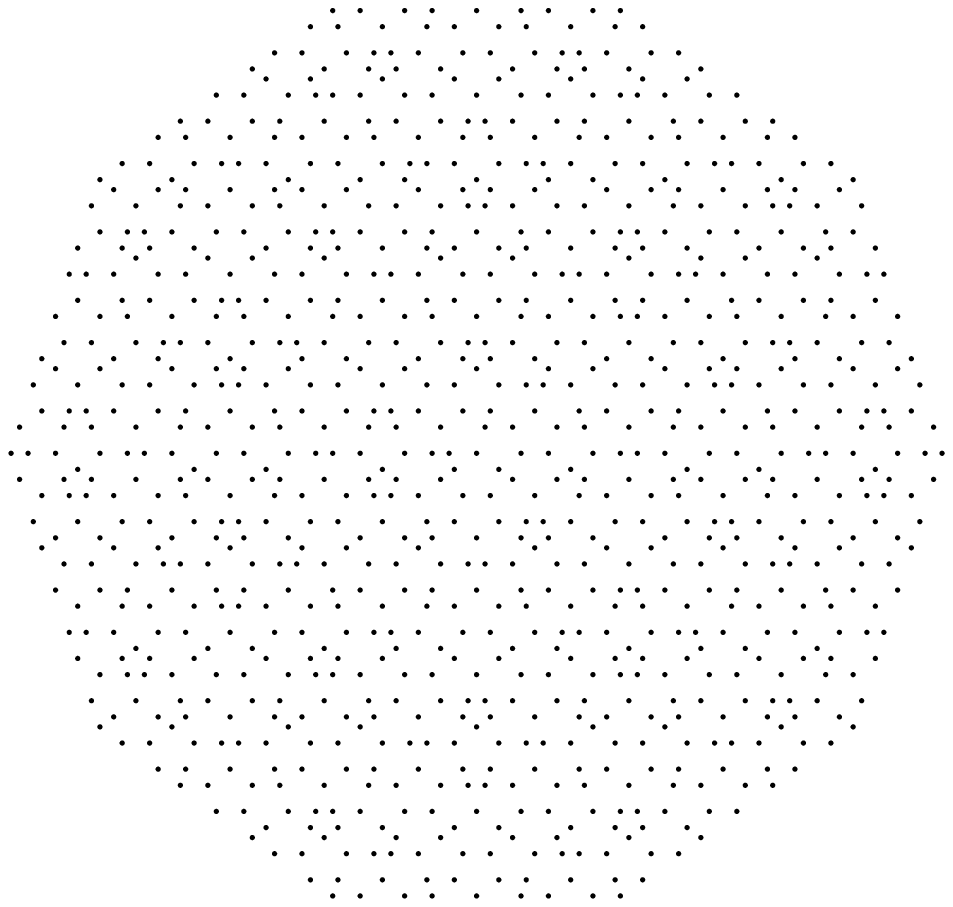}\hspace{0.1\textwidth}
\epsfysize=0.48\textwidth\epsfbox{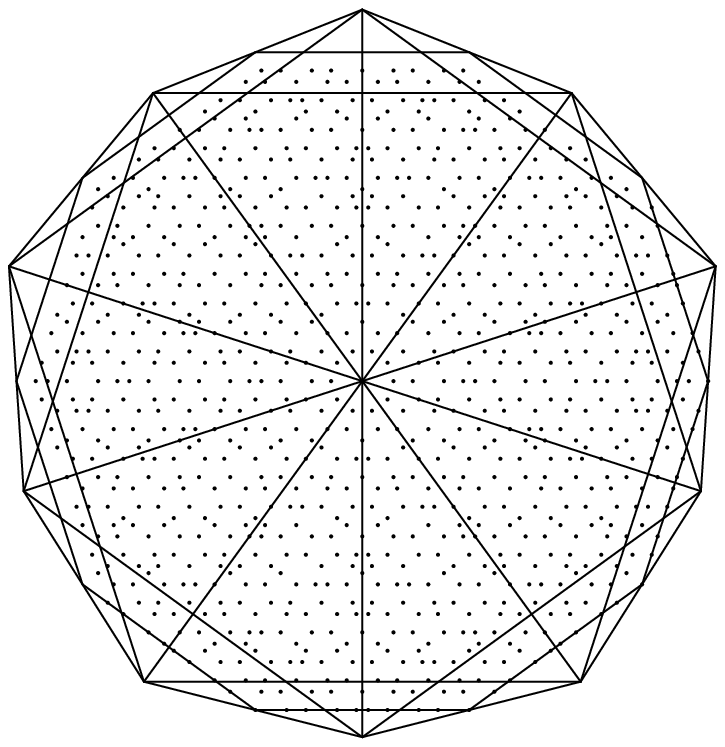}}
\caption{The central slice of the patch of $\varLambda_{\rm ico}^{\rm B}$ from Figure~\ref{fig:slices} (left) and its $.^{\star}$-image inside the (marked) decagon $(s+W)\cap H^{(\tau',0,1)}$ (right), both seen from perpendicular viewpoints.}
\label{fig:slices1}
\end{figure}

\begin{ex}
For illustrations of cyclotomic model sets with underlying
$\Z$-module $\Z[\zeta_5]$, see Figure~\ref{fig:slices1} on the left and Figure~\ref{fig:slices2}; {\it cf.} Proposition~\ref{slice} and Example~\ref{sliceex} below.
\end{ex}

\begin{lem}\label{icperp}
For $L\in\{\operatorname{Im}(\Ic),\Ic_0\}$, the following equations hold:
\begin{itemize}
\item[{\rm (a)}] 
$L\cap H^{(\tau,0,1)}=\Z[\tau](0,1,0)^{t}\oplus\Z[\tau]\tfrac{1}{2}(-1,-\tau',\tau)^{t}$.
\item[{\rm (b)}]
$
(L\cap H^{(\tau,0,1)})^{\star}=L^{\star}\,\,\cap\,\,H^{(\tau',0,1)}
$.
\end{itemize}
\end{lem}
\begin{proof}
Part (a) follows from Equations~\eqref{mbeqn} and~\eqref{mbeqn2}
together with the
relations $\operatorname{Im}(\Ic)=\tfrac{1}{2}\mathcal{M}_{\text{B}}$
and $\Ic_0=\tfrac{1}{2}\mathcal{M}_{\text{F}}$. Part (b)
follows from the identity
$((\tau,0,1)^t)^{\star}=(\tau',0,1)^t$.
\end{proof}

%\begin{coro}
%For all $\alpha\in\operatorname{Im}(\Ic)$, one has
%\begin{equation}\label{eq1}
%\begin{split}
%&\operatorname{Im}(\Ic)\,\,\cap\,\,\alpha+\left(\R(\tau,0,1)^{t}\right)^{\perp}\\ & 
%=\alpha+\big\langle(0,1,0)^{t},\tfrac{1}{2}(-1,-\tau',\tau)^{t}\big\rangle_{Z[\tau]}
%\end{split}
%\end{equation}
%and
%\begin{equation}\label{eq2}
%\begin{split}
%&{\operatorname{Im}(\Ic)}^{\star}\,\,\cap\,\,\alpha^{\star}+\left(\R(\tau',0,1)^{t}\right)^{\perp}\\ &=\alpha^{\star}+\big\langle(0,1,0)^{t},\tfrac{1}{2}(-1,-\tau,\tau')^{t}\big\rangle_{Z[\tau]}\,.
%\end{split}
%\end{equation}
%\end{coro}
%\begin{proof}
%Equation~\eqref{eq1} follows immediately from Lemma~\ref{icperp},
%while Equation~\eqref{eq2} follows from the identities
%$(\tau,0,1)^{\star}=(\tau',0,1)$, $(0,1,0)^{\star}=(0,1,0)$ and
%$(\tfrac{1}{2}(-1,-\tau',\tau))^{\star}=\tfrac{1}{2}(-1,-\tau,\tau')$ in conjunction with Equation~\eqref{eq1}.
%\end{proof}

\begin{defi}
We denote by
 $\Phi$ the $\R$-linear isomorphism  
$
\Phi\!:\,H^{(\tau,0,1)} \rightarrow \C
$, 
determined by $(0,1,0)^{t}\mapsto 1$ and $\tfrac{1}{2}(-1,-\tau',\tau)^{t}\mapsto\zeta_{5}$. Further, $\Phi^{\star}$ will denote the $\R$-linear isomorphism $
\Phi^{\star}\!:\,H^{(\tau',0,1)} \rightarrow \C
$, determined by $(0,1,0)^{t}\mapsto 1$ and $\tfrac{1}{2}(-1,-\tau,\tau')^{t}\mapsto \zeta_{5}^{3}$.
\end{defi}

\begin{lem}\label{remzm}
The maps $\Phi$
and $\Phi^{\star}$ are isometries of Euclidean vector
spaces, where $H^{(\tau,0,1)}$, $H^{(\tau',0,1)}$
and $\C$ are regarded as two-dimensional Euclidean vector spaces in the canonical way. Moreover, identifying   
$\C$ with the $xy$-plane in $\R^3$, $\Phi$ and $\Phi^{\star}$ extend uniquely to direct rigid motions of $\R^3$, {\it i.e.}, elements 
of the group $\operatorname{SO}(3,\R)$.  

\begin{proof}
The first assertion follows from the following identities:  
$$
\left\Arrowvert r(0,1,0)^{t}+s \tfrac{1}{2}(-1,-\tau',\tau)^{t}\right\Arrowvert=|r+s\,\zeta_{5}|
=\sqrt{r^2+s^2-rs\tau'}\,,
$$ $$\left\Arrowvert r(0,1,0)^{t}+s \tfrac{1}{2}(-1,-\tau,\tau')^{t} \right\Arrowvert=|r+s\,\zeta_{5}^{3}|=\sqrt{r^2+s^2-rs\tau}\,.$$
The additional statement is immediate.
%\begin{eqnarray*}
%\big\Arrowvert r(0,1,0)^{t}+s\tfrac{1}{2}(-1,-\tau',\tau)^{t} \big\Arrowvert&=&|r+s\,\zeta_{5}|\\
%&=&\sqrt{r^2+s^2-rs\tau'}
%\end{eqnarray*}
%and
%\begin{eqnarray*}
%\big\Arrowvert r(0,1,0)^{t}+s\tfrac{1}{2}(-1,-\tau,\tau')^{t}\big \Arrowvert&=&|r+s\,\zeta_{5}^{3}|\\
%&=&\sqrt{r^2+s^2-rs\tau}\,.
%\end{eqnarray*} 
\end{proof}
\end{lem}

\begin{lem}\label{res}
Let $L\in\{\operatorname{Im}(\Ic),\Ic_0\}$. Via
restriction, the maps $\Phi$ and $\Phi^{\star}$ 
induce isomorphisms of rank two $\Z[\tau]$-modules:
$$
L\cap H^{(\tau,0,1)}\stackrel{\Phi}{\longrightarrow}\Z[\zeta_5]\,,
$$
$$
L^{\star}\cap H^{(\tau',0,1)}\stackrel{\Phi^{\star}}{\longrightarrow}\Z[\zeta_5]\,.
$$

\end{lem}
\begin{proof}
This follows immediately from the definition of $\Phi$ and $\Phi^{\star}$ together with Lemma~\ref{icperp} and Remark~\ref{remzm2}.
\end{proof}

%\begin{coro}\label{zi}
%The set 
%\begin{eqnarray*}
%\mathcal{B}_{\operatorname{Im}(\Ic)}&:=&\big\{\tfrac{1}{2}(1,1,1)^t,\tfrac{\tau}{2}(1,1,1)^t,(0,1,0)^t,\tfrac{1}{2}(-1,-\tau',\tau)^t,\\
%&\hphantom{:=}&\hphantom{\big\{}\tfrac{1}{2}(\tau',-\tau,1)^t,\tfrac{1}{2}(-\tau',-\tau,-1)^t\big\}\,.
%\end{eqnarray*}
%is simultaneously a $\Z$-basis of $\operatorname{Im}(\Ic)$ and a $\Q$-basis of $\langle\operatorname{Im}(\Ic)\rangle_{\Q}=\Q\operatorname{Im}(\Ic)=(\mathbbm{Q}(\tau))^3$. 
%\end{coro}
%\begin{proof}
%First, observe the equalities $$\Phi((0,1,0)^t)=1\,,\,\,
%\Phi(\tfrac{1}{2}(-1,-\tau',\tau)^t)=\zeta_5\,,\,\,
%\Phi(\tfrac{1}{2}(\tau',-\tau,1)^t)=\zeta_5^2$$ and
%$\Phi(\tfrac{1}{2}(-\tau',-\tau,-1)^t)=\zeta_5^3$. Then, the fact
%that $\mathcal{B}_{\operatorname{Im}(\Ic)}$ is simultaneously
%a $\Z$-basis of $\operatorname{Im}(\Ic)$ and a $\Q$-basis of
% $\langle\operatorname{Im}(\Ic)\rangle_{\Q}=\Q\operatorname{Im}(\Ic)$ follows immediately from Lemma~\ref{res} and Equation~\eqref{mbeqn} in conjunction with the
%fact that the set $\{1,\zeta_5,\zeta_5^2,\zeta_5^3\}$ is simultaneously
%a $\Z$-basis
%of the $\Z$-module $\Z[\zeta_5]$ and a $\Q$-basis of $\mathbbm{Q}(\zeta_5)=\langle\Z[\zeta_5]\rangle_{\Q}=\Q\Z[\zeta_5]$. The equality
%$(\mathbbm{Q}(\tau))^3=\langle\operatorname{Im}(\Ic)\rangle_{\Q}$
% then follows immediately from the obvious inclusion
%$\langle\operatorname{Im}(\Ic)\rangle_{\Q}\subset (\mathbbm{Q}(\tau))^3$
%in conjunction with the fact that both sets are $6$-dimensional vector
%spaces over $\Q$.
%\end{proof}

\begin{prop}\label{slice}
Let 
$\varLambda$ be a generic icosahedral model set with underlying
$\Z$-module $L$, say
$\varLambda=\varLambda_{{\rm ico}}(t,W)$. Then, for every
$\lambda\in\varLambda$, one has the identity 
$$
\Phi\big((\varLambda\,\,\cap\,\,(\lambda+H^{(\tau,0,1)}))-\lambda\big) \,\,=\,\, \big\{z\in\Z[\zeta_5]\,\big|\,z^{\star_{5}}\in W_{\lambda}\big\}\,,
$$
where $.^{\star_{5}}$ is the Galois automorphism of $\mathbbm{Q}(\zeta_5)/\Q$, defined by $\zeta_{5}\mapsto
\zeta_{5}^{3}$ and 
$$
W_{\lambda}:=\Phi^{\star}\big((W\,\,\cap\,\,((\lambda-t)^{\star}+H^{(\tau',0,1)}))-(\lambda-t)^{\star}\big)\,.
$$
Thus, the sets of the form
\begin{eqnarray}\label{slice2}
\Phi\big((\varLambda\,\,\cap\,\,(\lambda+H^{(\tau,0,1)}))-\lambda\big)\,,
\end{eqnarray}
where $\lambda\in\varLambda$, are cyclotomic model sets with underlying
$\Z$-module $\Z[\zeta_5]$.
\end{prop}
\begin{proof}
First, consider $\Phi(\mu)$, where $\mu\in (\varLambda\cap
    (\lambda+H^{(\tau,0,1)}))-\lambda$. It follows that
    $\mu\in L\cap H^{(\tau,0,1)}$ and
    $(\mu+(\lambda-t))^{\star}=\mu^{\star}+(\lambda-t)^{\star}\in W$. Lemma~\ref{res} implies that $\Phi(\mu)\in\Z[\zeta_5]$, say
    $\Phi(\mu)=\alpha+\beta\zeta_5$ for suitable
    $\alpha,\beta\in\Z[\tau]$. One has
$$
\Phi(\mu)^{\star_5}\,=\,\alpha'+\beta'\zeta_5^{3}\,=\,\Phi^{\star}(\mu^{\star})\,\,\in\,\, W_{\lambda}\,.
$$
Conversely, suppose that $z\in\Z[\zeta_5]$ satisfies $z^{\star_5}\in
W_{\lambda}$. Then, there are suitable
    $\alpha,\beta\in \Z[\tau]$ such that $z=\alpha+\beta\zeta_5$ and,
    consequently, 
    $z^{\star_5}=\alpha'+\beta'\zeta_5^{3}\in W_{\lambda}$. By
    definition of $W_{\lambda}$, one has
    $z^{\star_5}=\Phi^{\star}(\mu)$, where $\mu\in
    H^{(\tau',0,1)}$ satisfies $\mu+(\lambda-t)^{\star}\in
    W$. Clearly, there exist $r,s\in\R$ such that $\mu=r(0,1,0)^t+s
    \tfrac{1}{2}(-1,-\tau,\tau')^t$, whence
    $\Phi^{\star}(\mu)=r+s\zeta_5^{3}$. The linear independence
    of $1$ and $\zeta_5^{3}$ over $\R$ now implies that $r=\alpha$ and
    $s=\beta$, so that
    $\mu\in L^{\star}$. Moreover, one can verify that one has $\mu^{-\star}\in (\varLambda\cap
    (\lambda+H^{(\tau,0,1)}))-\lambda$ and $\Phi(\mu^{-\star})=\alpha+\beta\zeta_5=z$. This proves the claimed identity. The assertion is now immediate.
\end{proof}

\begin{ex}\label{sliceex}
For an illustration of the content of Proposition~\ref{slice} in case of the icosahedral model set $\varLambda_{\rm ico}^{\rm B}$ from Example~\ref{lambdaico}, see Figures~\ref{fig:slices1} and~\ref{fig:slices2}.
\end{ex}

\begin{figure}
\centerline{\epsfysize=0.48\textwidth\epsfbox{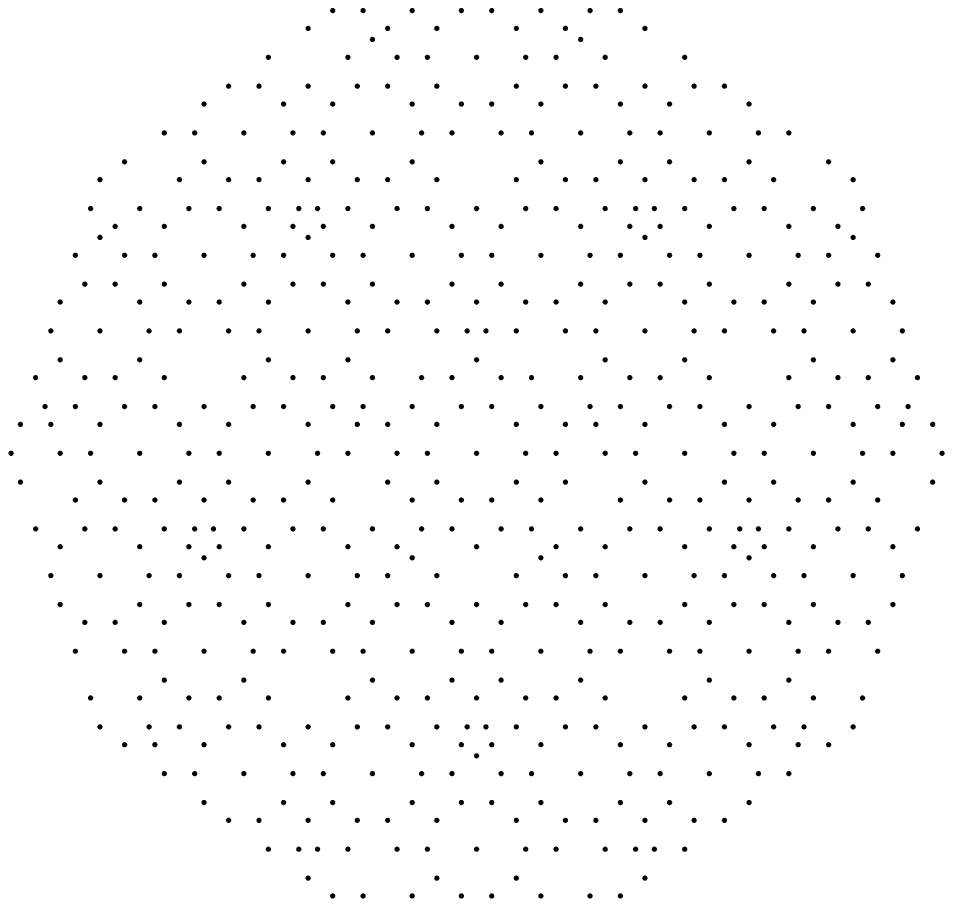}\hspace{0.01\textwidth}
\epsfysize=0.48\textwidth\epsfbox{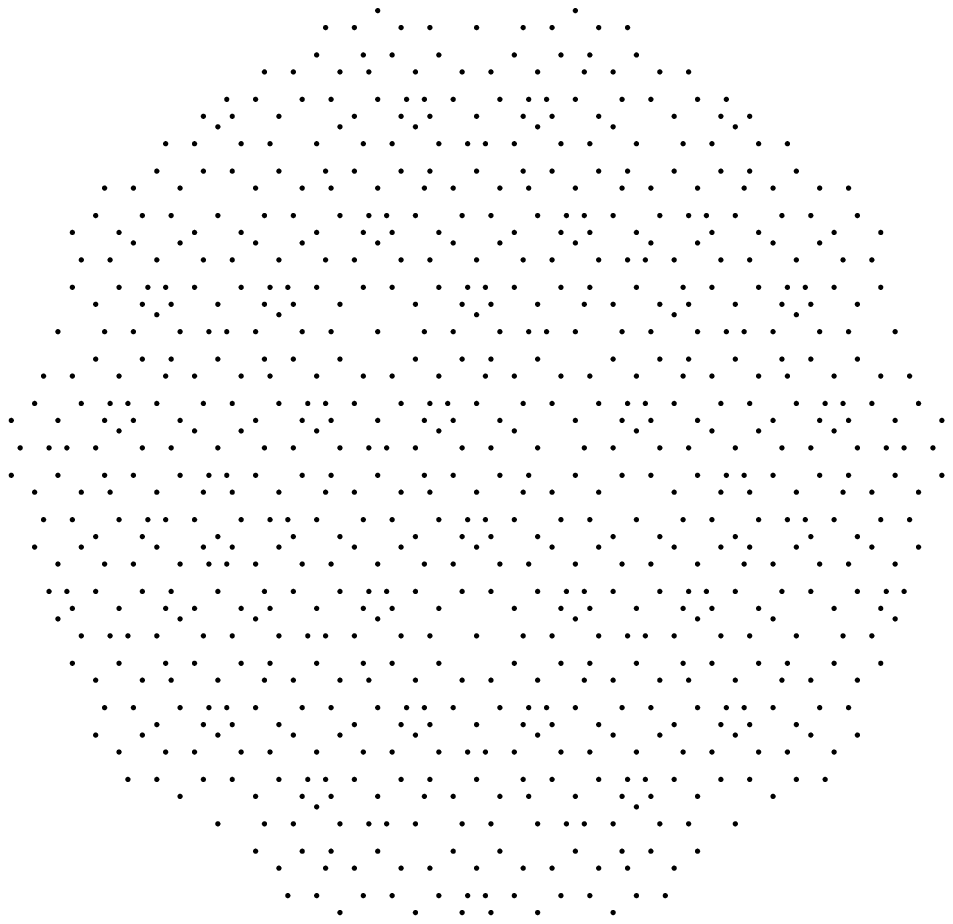}}
\caption{Another two slices of the patch of $\varLambda_{\rm ico}^{\rm B}$ from Figure~\ref{fig:slices}.}
\label{fig:slices2}
\end{figure}

%\begin{figure}
%\centerline{\epsfysize=0.42\textwidth\epsfbox{vp4a.eps}\hspace{0.03\textwidth}
%\epsfysize=0.42\textwidth\epsfbox{vp5a.eps}}
%\caption{A few slices of a patch of the icosahedral model set $\varLambda_{\rm ico}$.}
%\label{fig:slices3}
%\end{figure}

\subsection{The translation module of icosahedral model sets}

In order to shed some light on the set of $\varLambda$-directions of
  an icosahedral model set $\varLambda$ with underlying
$\Z$-module $L$, we first have to establish a relation between icosahedral model sets and their
  underlying $\Z$-modules. We denote by $m_{\tau}$ the $\Z[\tau]$-module endomorphism of
$\Q(\tau)^3$, given by multiplication by $\tau$, {\it i.e.}, $\alpha\mapsto
\tau\alpha$. Furthermore, we denote by ${m_{\tau}}^{\star}$ the $\Z[\tau]$-module endomorphism of
$(\Q(\tau)^3)^{\star}$, given by $\alpha^{\star}\mapsto
(\tau\alpha)^{\star}$.

\begin{lem}\label{r2ico}
The map
 ${m_{\tau}}^{\star}$ is contractive with contraction constant
 $1/\tau\in(0,1)$, {\em i.e.}, the equality $\Arrowvert
 {m_{\tau}}^{\star}(\alpha^{\star})\Arrowvert = (1/\tau)\, \Arrowvert
 \alpha^{\star}\Arrowvert$ holds for all $\alpha\in \Q(\tau)^3$.
\end{lem}
\begin{proof}
For $\alpha\in \Q(\tau)^3$, observe that 
$
\Arrowvert {m_{\tau}}^{\star}(\alpha^{\star})\Arrowvert=\Arrowvert
(\tau \alpha)^{\star}\Arrowvert = \Arrowvert
\tau' \alpha^{\star}\Arrowvert = (1/\tau) \, \Arrowvert
\alpha^{\star}\Arrowvert$. 
\end{proof}

\begin{lem}\label{dilateico}
Let $\varLambda$ be an icosahedral model
set with underlying
$\Z$-module $L$, say $\varLambda=\varLambda_{{\rm ico}}(t,W)$. Then, for any $F\in\mathcal{F}(L)$, there is an expansive homothety $h\!:\, \R^3 \rightarrow
\R^3$ such that $h(F)\subset \varLambda$. 
\end{lem}
\begin{proof}
From
$\operatorname{int}(W)\neq \varnothing$ and the denseness of $L^{\star}$ in
$\R^3$, one gets the existence of a
suitable $\alpha_{0}\in L$ with ${\alpha_{0}}^{\star}\in\, \operatorname{int}(W)$. Consider the open neighbourhood $V:= \,\operatorname{int}(W) - \, \alpha_{0}{}^{\star}$ of $0$ in
$\R^3$. Since the map ${m_{\tau}}^{\star}$ is
contractive by Lemma~\ref{r2ico} (in the sense which was made precise in that lemma),
 the existence of a suitable $k\in\mathbbm{N}$ is implied such that
$({m_{\tau}}^{\star})^{k}(F^{\star})\subset
V$. Hence, one has $\{(\tau^{k} \alpha +
\alpha_{0})^{\star}\, |\, \alpha\in F\}\subset
\operatorname{int}(W)\subset W$ and, further, $h(F)\subset \varLambda$, where $h\!:\, \R^3
\rightarrow \R^3$ is the expansive homothety given by $x \mapsto
\tau^{k} x + (\alpha_{0}+t)$. 
\end{proof}

As an easy application of Lemma~\ref{dilateico}, one obtains the
following result on the set of $\varLambda$-directions for icosahedral model
sets $\varLambda$.

\begin{prop}\label{direc}
Let $\varLambda$ be an icosahedral model
set with underlying
$\Z$-module $L$. Then, the set of $\varLambda$-directions is
precisely the set of $L$-directions.
\end{prop}
\begin{proof}
Since one has $\varLambda - \varLambda\subset L$, every $\varLambda$-direction
is an $L$-direction. For the converse, let $u\in
\mathbb{S}^{2}$ be an $L$-direction, say parallel to
 $\alpha\in L\setminus \{0\}$. By Lemma~\ref{dilateico}, there is a homothety $h\!:\, \R^3 \rightarrow
\R^3$ such that $h(\{0,\alpha\})\subset \varLambda$. It follows that
$h(\alpha)-h(0)\in (\varLambda -
\varLambda)\setminus\{0\}$. Since $h(\alpha)-h(0)$ is parallel to
$\alpha$, the assertion follows. 
\end{proof}

\section{Complexity}\label{secset}

In the practice of quantitative HRTEM, the determination of the rotational orientation of a
quasicrystalline probe in an electron microscope can 
rather easily be achieved in the diffraction mode. This is due to the icosahedral symmetry of genuine
 icosahedral quasicrystals. However, the $X$-ray images  taken
in the high-resolution mode do not allow us to locate the examined sets. Therefore, as already pointed out in~\cite{BG2}, in order to
prove practically relevant and rigorous results, one has to deal
with the {\em non-anchored case} of the whole {\em local indistinguishability class} (or LI-class, for short) $\operatorname{LI}(\varLambda)$ of a regular, generic icosahedral model set
$\varLambda$, rather than dealing with the {\em anchored case} of a single fixed icosahedral model set $\varLambda$; recall Remark~\ref{propcm} for the equivalence relation given by local indistinguishability and compare also~\cite{GL}.

\begin{rem}
In the crystallographic case of a lattice $L$ in $\R^3$, the LI-class of $L$ consists of all translates of $L$ in $\R^3$, {\it i.e.}, one has $\operatorname{LI}(L)=\{t+L\,|\,t\in\R^3\}$. In particular, $\operatorname{LI}(L)$ simply consists of one translation class. The entire LI-class $\operatorname{LI}(\varLambda_{{\rm ico}}(t,W))$ of a regular, generic icosahedral model set $\varLambda_{{\rm ico}}(t,W)$
can be shown to consist of all generic icosahedral model sets of the form $\varLambda_{{\rm ico}}(t,s+W)$ and all patterns obtained as limits of sequences of generic icosahedral model sets of the form $\varLambda_{{\rm ico}}(t,s+W)$ in the local
topology (LT). Here, two patterns are $\varepsilon$-close if, after a
translation by a distance of at most $\varepsilon$, they agree on a ball of
radius $1/\varepsilon$ around the origin; see~\cite{B,Schl}. Each such limit is then a subset of some
$\varLambda_{{\rm ico}}(t,s+W)$, but $s$ might not be in a generic
position. Note that the LI-class $\operatorname{LI}(\varLambda)$ of an icosahedral model set $\varLambda$ contains {\em uncountably many} (more precisely, $2^{\aleph_0}$) translation classes; {\it cf.}~\cite{B} and references therein.
\end{rem}

In view of the complication described above, we must make sure that we
deal with finite subsets of {\em generic} icosahedral model sets of the form $\varLambda_{{\rm ico}}(t,s+W)$, {\it
i.e.}, subsets whose $.^{\star}$-image lies in the {\em interior} of
the window. This
restriction to the generic case is the proper analogue of the
restriction to {\em perfect} lattices and their translates in the
crystallographic case. Analogous to the lattice case~\cite{GG2,GGP} and the case of cyclotomic model sets~\cite{BG2}, the main algorithmic problems of the discrete tomography of icosahedral model sets look as follows.

\begin{defi}[Consistency, Reconstruction, and Uniqueness
  Problem]\label{m-def:consrecprob}
Let $L=\operatorname{Im}(\Ic)$ (resp., $L=\Ic_0$), let $W\subset \mathbbm{R}^3$ be a window and let
  $u_1,\dots,u_m\in\mathbb{S}^2$ be $m\geq 2$ pairwise non-parallel
  $L$-directions. The corresponding consistency,
  reconstruction and uniqueness problems are defined as follows.

  \begin{quote}
    {\sc Consistency}. \\
    Given functions $p_{u_{j}} : \mathcal{L}_{u_{j}}^3 \rightarrow
    \mathbbm{N}_{0}$, $j\in\{1,\dots,m\}$, whose supports are finite
    and satisfy $\operatorname{supp}(p_{u_{j}})\subset
    \mathcal{L}^{L}_{u_{j}}$, decide whether there is a
    finite set $F$ which is contained in an element of 
    $\mathcal{I}^{\rm B}_g(W)$ (resp., $\mathcal{I}^{\rm F}_g(W)$) and satisfies
    $X_{u_{j}}F=p_{u_{j}}$, $j\in\{1,\dots,m\}$.
\end{quote}

\begin{quote}
  {\sc Reconstruction}. \\
  Given functions $p_{u_{j}} : \mathcal{L}_{u_{j}}^3 \rightarrow
  \mathbbm{N}_{0}$, $j\in\{1,\dots,m\}$, whose supports are finite and
  satisfy $\operatorname{supp}(p_{u_{j}})\subset
  \mathcal{L}^{L}_{u_{j}}$, decide whether there exists
  a finite subset $F$ of an element of 
  $\mathcal{I}^{\rm B}_g(W)$ (resp., $\mathcal{I}^{\rm F}_g(W)$) that satisfies
  $X_{u_{j}}F=p_{u_{j}}$, $j\in\{1,\dots,m\}$, and, if so, construct
  one such $F$.
\end{quote}

\begin{quote}
  {\sc Uniqueness}. \\
  Given a finite subset $F$ of an element of
  $\mathcal{I}^{\rm B}_g(W)$ (resp., $\mathcal{I}^{\rm F}_g(W)$), decide whether
  there is a different finite set $F'$ that is also a subset of an  element of 
  $\mathcal{I}^{\rm B}_g(W)$ (resp., $\mathcal{I}^{\rm F}_g(W)$) and satisfies
  $X_{u_{j}}F=X_{u_{j}}F'$, $j\in\{1,\dots,m\}$.
\end{quote}

\end{defi}

One has the following tractability result, which was proved for the
case of B-type icosahedral model sets by combining the results from
Section~\ref{cycsec} with those presented in~\cite{BG2}; {\it
  cf.}~\cite[Theorem 3.33]{H2} for the details. The proof for the
F-type case is similar and we prefer to omit the straightforward details
here. Below, for $L\in\{\operatorname{Im}(\Ic),\Ic_0\}$, the $L$-directions in $\mathbb{S}^2\cap
H^{(\tau,0,1)}$ will be called \emph{$L^{(\tau,0,1)}$-directions}. By
Lemma~\ref{icperp}(a), the set of
$\operatorname{Im}(\Ic)^{(\tau,0,1)}$-directions and the set of
$\Ic_0^{(\tau,0,1)}$-directions coincide.

\begin{theorem}\label{th2}
Let $L\in\{\operatorname{Im}(\Ic),\Ic_0\}$. When restricted to two $L^{(\tau,0,1)}$-directions and polyhedral
  windows, the problems {\sc Consistency}, {\sc Reconstruction} and {\sc
  Uniqueness} as defined in Definition~\ref{m-def:consrecprob} 
  can be solved in polynomial time in the real RAM-model of computation.
\end{theorem}

\begin{rem}
For a detailed analysis of the complexities of the above algorithmic
problems in the B-type case, we refer the reader to~\cite[Chapter 3]{H2}. Note 
that even in the anchored planar lattice case $\mathbbm{Z}^2$ the corresponding problems {\sc
  Consistency}, {\sc Reconstruction} and {\sc Uniqueness} 
are $\mathbbm{NP}$-hard for three or more $\Z^2$-directions;
{\it cf.}~\cite{GG2,GGP}.
\end{rem}

\section{Uniqueness}

\subsection{Simple results on determination of finite subsets of icosahedral model sets}\label{simple}

In this section, we present some uniqueness results which only deal
with the {\em anchored case} of determining finite subsets of a fixed
icosahedral model set $\varLambda$ by $X$-rays in
{\em arbitrary} $\varLambda$-directions; {\it cf.}
Proposition~\ref{direc}. As already explained in Section~\ref{intro},
$X$-rays in non-$\varLambda$-directions are meaningless in practice. Without the restriction to $\varLambda$-directions, the finite subsets of a fixed icosahedral 
model set $\varLambda$ can be determined by one $X$-ray. In fact, any $X$-ray in a non-$\varLambda$-direction is suitable
for this purpose, since any line in $3$-space in a
non-$\varLambda$-direction passes through at most one point of
$\varLambda$.      
The next result represents a fundamental source of difficulties
in discrete tomography. There exist several versions;
compare~\cite[Theorem 4.3.1]{HK}, \cite[Lemma 2.3.2]{G},
\cite[Proposition 4.3]{BH}, \cite[Proposition 2.3 and Remark 2.4]{H2} and
\cite[Proposition 8]{H}. 

\begin{prop}\label{source}
Let $\varLambda$ be an icosahedral model
set with underlying
$\Z$-module $L$, say $\varLambda=\varLambda_{{\rm ico}}(t,W)$. Further, let $U\subset \mathbb{S}^{2}$ be an arbitrary, but fixed finite set of pairwise non-parallel $L$-directions. Then, $\mathcal{F}(\varLambda)$ is not determined by the $X$-rays in the directions of $U$.
\end{prop}
\begin{proof}
We argue by induction on $\operatorname{card}(U)$. The case
$\operatorname{card}(U)=0$ means $U=\varnothing$ and is obvious. Suppose the assertion to be
true whenever $\operatorname{card}(U)=k\in \mathbbm{N}_{0}$ and let
$\operatorname{card}(U)=k+1$. By induction hypothesis, there are
different elements $F$ and $F'$ of $\mathcal{F}(\varLambda)$
with the same $X$-rays in the directions of $U'$, where $U'\subset U$
satisfies $\operatorname{card}(U')=k$. Let $u$ be the remaining
direction of $U$. Choose a non-zero element $\alpha\in L$ parallel to
$u$ such that $\alpha+(F\cup F')$ and $F\cup F'$ are disjoint. Then,
$F'':= (F\cup (\alpha+F'))-t$ and
$F''':= (F'\cup (\alpha+F))-t$ are different
elements of $\mathcal{F}(L)$ with the same $X$-rays in the
directions of $U$. By Lemma~\ref{dilateico}, there is a homothety $h\!:\,
\R^3 \rightarrow \R^3$ such that $h(F''\cup F''')=h(F'')\cup
h(F''')\subset \varLambda$. It follows that $h(F'')$ and
$h(F''')$ are different elements of $\mathcal{F}(\varLambda)$
with the same $X$-rays in the directions of $U$; {\it cf.}~Lemma~\ref{homotu}. 
\end{proof}

\begin{rem}\label{source2}
An analysis of the proof of Proposition~\ref{source} shows that, for
any finite set $U\subset \mathbb{S}^{2}$ of $k$ pairwise non-parallel
$L$-directions, there are disjoint elements  $F$ and $F'$ of
$\mathcal{F}(\varLambda)$ with
$\operatorname{card}(F)=\operatorname{card}(F')=2^{(k-1)}$ and with the same
$X$-rays in the directions of $U$. Consider any convex subset $C$ of 
$\R^3$ which contains $F$ and $F'$ from above. Then, the subsets
$F_1:=(C\cap\varLambda)\setminus F$ and
$F_2:=(C\cap\varLambda)\setminus F'$ of
$\mathcal{F}(\varLambda)$ also have the same
$X$-rays in the directions of $U$. Whereas the points in $F$ and $F'$
are widely dispersed over a region, those in $F_1$ and $F_2$ are
contiguous in a way similar to atoms in a quasicrystal;
compare~\cite[Remark 4.3.2]{GG2} and~\cite[Remark 2.4 and Figure 2.1]{H} (see also~\cite[Remark 32 and Figure 5]{H}). 
\end{rem}

Originally, the proof of the following result is due to R\'enyi; {\it cf.}~\cite{Re} and compare~\cite[Theorem
4.3.3]{HK}.

\begin{prop}\label{m+1}
Let $\varLambda$ be an icosahedral model
set with underlying
$\Z$-module $L$. Further, let $U\subset \mathbb{S}^{2}$ be any set of $k+1$ pairwise non-parallel $L$-directions, where $k\in \mathbbm{N}_{0}$. Then, $\mathcal{F}_{\leq k}(\varLambda)$ is determined by the $X$-rays in the directions of $U$. Moreover, for all $F\in\mathcal{F}_{\leq
  k}(\varLambda)$, one has $G^{F}_{U}=F$.
\end{prop}
\begin{proof}
Let $F,F'\in \mathcal{F}_{\leq k}(\varLambda)$ have the same $X$-rays in the directions of $U$. Then, one has $\operatorname{card}(F)=\operatorname{card}(F')$ by Lemma~\ref{cardinality}(a) and $
F,F'\subset G_{F}^{U}$ by Lemma~\ref{fgrid}. But we have $G_{F}^{U}=F$ since the existence of a point in $G_{F}^{U}\setminus F$ implies the existence of at least $\operatorname{card}(U)\geq k+1$ points in $F$, a contradiction. It follows that $F=F'$.
\end{proof}

\begin{rem}\label{exdodecagon}
In particular, the additional statement of Proposition~\ref{m+1}
demonstrates that, for a fixed icosahedral model set $\varLambda$ with underlying
$\Z$-module $L$, the unique reconstruction of sets $F\in\mathcal{F}_{\leq k}(\varLambda)$ from their $X$-rays in arbitrary sets of $k+1$ pairwise non-parallel $L$-directions $U\subset \mathbb{S}^{2}$ merely amounts to compute the grids $G_{F}^{U}$. Let
$\varLambda$ be an icosahedral model set with underlying
$\Z$-module $L$. Remark~\ref{source2} and Proposition~\ref{m+1} show
that $\mathcal{F}_{\leq k}(\varLambda)$ can be determined by
the $X$-rays in
any set of $k+1$ pairwise non-parallel $L$-directions but not by
$1+\lfloor\log_{2}k\rfloor$ pairwise non-parallel $X$-rays in
$L$-directions. However, in practice, one is
interested in the determination of finite sets by $X$-rays in a small
number of directions since after about $3$ to $5$ images taken by
HRTEM, the object may be damaged or even destroyed by the radiation
energy. Observing that the typical atomic structures to be determined
comprise about $10^6$ to $10^9$ atoms, one realizes that the last
result is not practical at all.
\end{rem}

The following result was proved in~\cite[Theorem 2.8(a)]{H2}; see also~\cite[Theorem 13(a)]{H}.

\begin{prop}\label{flc}
Let $d\geq 2$, let $R>0$, and let $\varLambda\subset\R^d$ be a Delone set of finite local complexity. Then, the set $\mathcal{D}_{<R}(\varLambda)$ is determined by two $X$-rays in $\varLambda$-directions.
\end{prop}

Since icosahedral model sets $\varLambda\subset\R^3$ are Delone sets of finite local complexity ({\it cf.} Remark~\ref{propcm}), the following corollary follows immediately from Proposition~\ref{flc} in conjunction
with Proposition~\ref{direc}.

\begin{coro}\label{bounded}
Let
$\varLambda$ be an icosahedral model set with underlying
$\Z$-module $L$ and let $R>0$. Then, the set $\mathcal{D}_{<R}(\varLambda)$ is determined by two $X$-rays in $L$-directions.
\end{coro}

\begin{rem}
Although looking promising at first sight, Corollary~\ref{bounded} is of limited use in practice
because, in general, one cannot guarantee that all the directions
which are used yield densely occupied lines in icosahedral model sets.
\end{rem}

\subsection{Determination of convex subsets of icosahedral model sets}\label{convdet}

\begin{rem}\label{dirrem3}
Proposition~\ref{direc} shows that, for all icosahedral model sets
$\varLambda$ with underlying
$\Z$-module $L$, the set of $L^{(\tau,0,1)}$-directions is precisely the set of $\varLambda$-directions in $\mathbb{S}^2\cap
H^{(\tau,0,1)}$. Further, by Lemmas~\ref{remzm} and~\ref{res}, the set of $L^{(\tau,0,1)}$-directions maps under $\Phi$ bijectively onto the set of $\Z[\zeta_5]$-directions.
\end{rem}

The following property is evident.

\begin{lem}\label{trick7}
Let $L\in\{\operatorname{Im}(\Ic),\Ic_0\}$, let $U\subset \mathbb{S}^2$ be a finite set of $L^{(\tau,0,1)}$-directions, and let $F,F'\in\mathcal{F}(t+H^{(\tau,0,1)})$, where $t\in\R^3$. If $F$ and $F'$ have the same $X$-rays in the directions of $U$, then $\Phi(F-t)$ and $\Phi(F'-t)$ have the same $X$-rays in the directions of $\Phi(U)\subset \mathbb{S}^1$.
\end{lem}

The following fundamental result follows immediately from~\cite[Theorem 2.54]{H2}; see also~\cite[Theorem 15]{H}.

\begin{theorem}\label{convcyc}
The following assertions hold: 
\begin{itemize}
\item[{\rm (a)}]
There is a set $U\subset\mathbb{S}^1$ of four pairwise non-parallel
$\Z[\zeta_5]$-directions such that, for all cyclotomic model sets
$\varLambda_{\rm cyc}$ with underlying
$\Z$-module $\Z[\zeta_5]$, the set $\mathcal{C}(\varLambda_{\rm cyc})$ is determined by the $X$-rays in the directions of $U$.
\item[{\rm (b)}]
For all cyclotomic model sets $\varLambda_{\rm cyc}$ with underlying
$\Z$-module $\Z[\zeta_5]$ and all sets  $U\subset\mathbb{S}^1$ of three or less pairwise non-parallel $\Z[\zeta_5]$-directions, the set $\mathcal{C}(\varLambda_{\rm cyc})$ is not determined by the $X$-rays in the directions of $U$.
\end{itemize}
\end{theorem}

We are now able to prove the main result of this text by applying the results of~\cite{H2,H} on the determination of convex subsets of cyclotomic model sets with underlying
$\Z$-module $\Z[\zeta_5]$ to the various images $\Phi((\varLambda\cap (\lambda +H^{(\tau,0,1)}))-\lambda)$, where $\varLambda$ is an icosahedral model set and $\lambda\in\varLambda$.

\begin{rem}\label{trick8}
Note that, for a convex subset
$C$ of an icosahedral model set $\varLambda$ and an element $\lambda\in\varLambda$, the intersection $C\cap (\lambda+H^{(\tau,0,1)})$ is a convex subset of the slice $\varLambda\cap (\lambda+H^{(\tau,0,1)})$ of $\varLambda$. Hence, $\Phi((C\cap (\lambda+H^{(\tau,0,1)}))-\lambda)$ is a convex subset of $\Phi((\varLambda\cap (\lambda+H^{(\tau,0,1)}))-\lambda)$.
\end{rem}

The following fundamental result deals with the anchored case.

\begin{theorem}\label{main1ico}
Let $L\in\{\operatorname{Im}(\Ic),\Ic_0\}$. The following assertions hold:
\begin{itemize}
\item[{\rm (a)}]
There is a set $U\subset \mathbb{S}^2$ of four $L^{(\tau,0,1)}$-directions such that, for all generic icosahedral model sets
$\varLambda$ with underlying
$\Z$-module $L$, the set
$\mathcal{C}(\varLambda)$ is determined by the $X$-rays in
the directions of $U$.
\item[{\rm (b)}]
For all generic icosahedral model sets
$\varLambda$ with underlying
$\Z$-module $L$ and all sets $U\subset \mathbb{S}^2$  of three or less pairwise
non-parallel $L^{(\tau,0,1)}$-directions, the set
$\mathcal{C}(\varLambda)$ is not determined by the $X$-rays in the
directions of $U$. 
\end{itemize}
 \end{theorem}
\begin{proof}
For part (a), let $U'\subset\mathbb{S}^1$ be a set of four pairwise non-parallel $\Z[\zeta_5]$-directions with the property that, for all cyclotomic model sets $\varLambda_{\rm cyc}$ with underlying
$\Z$-module $\Z[\zeta_5]$, the set $\mathcal{C}(\varLambda_{\rm cyc})$ is determined by the $X$-rays in the directions of $U'$. Such a set $U'$ exists by Theorem~\ref{convcyc}(a). We claim that, for all generic icosahedral model sets
$\varLambda$ with underlying
$\Z$-module $L$, the set
$\mathcal{C}(\varLambda)$ is determined by the $X$-rays in
the directions of $U:=\Phi^{-1}(U')\subset\mathbb{S}^2$. {\it Cf.} Remark~\ref{dirrem3} for the fact that $U$ consists only of $L^{(\tau,0,1)}$-directions. Assume the existence of two different elements, say $C$ and $C'$, of $\mathcal{C}(\varLambda)$ having the same $X$-rays in the directions of $U$. Hence, there is an element $\lambda\in\varLambda$ such that $C\cap (\lambda+H^{(\tau,0,1)})$ and $C'\cap (\lambda+H^{(\tau,0,1)})$ are different convex subsets of the slice $\varLambda\cap (\lambda+H^{(\tau,0,1)})$ with the same $X$-rays in the directions of $U$. By Lemma~\ref{trick7} and Remark~\ref{trick8}, it follows that $\Phi((C\cap (\lambda+H^{(\tau,0,1)}))-\lambda)$ and $\Phi((C\cap (\lambda+H^{(\tau,0,1)}))-\lambda)$ are different convex subsets of $\Phi((\varLambda\cap (\lambda+H^{(\tau,0,1)}))-\lambda)$ having the same $X$-rays in the $\Z[\zeta_5]$-directions of $U'$. Since the set $\Phi((\varLambda\cap (\lambda+H^{(\tau,0,1)}))-\lambda)$ is a cyclotomic model set with underlying
$\Z$-module $\Z[\zeta_5]$ by Proposition~\ref{slice}, this is a contradiction. 

For assertion (b), let $U\subset \mathbb{S}^2$ be a set of three or less pairwise
non-parallel $L^{(\tau,0,1)}$-directions and let $\varLambda$ be a
generic icosahedral model set with underlying
$\Z$-module $L$. Consider a slice $\varLambda\cap (\lambda+H^{(\tau,0,1)})$ of $\varLambda$, $\lambda\in\varLambda$, together with the cyclotomic model set $\Phi((\varLambda\cap (\lambda+H^{(\tau,0,1)}))-\lambda)$ with underlying
$\Z$-module $\Z[\zeta_5]$; {\it cf.} Proposition~\ref{slice}. By Theorem~\ref{convcyc}(b), there are two different convex subsets, say $C$ and $C'$, of $\Phi((\varLambda\cap (\lambda+H^{(\tau,0,1)}))-\lambda)$ with the same $X$-rays in the $\Z[\zeta_5]$-directions of $U':=\Phi(U)\subset\mathbb{S}^1$; {\it cf.} Remark~\ref{dirrem3}. It follows that $\Phi^{-1}(C)+\lambda$ and $\Phi^{-1}(C')+\lambda$ are different convex subsets of (the slice $\varLambda\cap (\lambda+H^{(\tau,0,1)})$ of) $\varLambda$ with the same $X$-rays in the $L^{(\tau,0,1)}$-directions of $U$.  
\end{proof}

\begin{rem}\label{nonanch2}
The proof of Theorem~\ref{main1ico} shows that the result extends
to the set of subsets $C$ of generic icosahedral model sets $\varLambda$ that are only $H^{(\tau,0,1)}${\em -convex}, the latter
meaning that, for all
$\lambda\in \varLambda$, the sets $C\cap (\lambda+H^{(\tau,0,1)})$ are convex subsets of the slices $\varLambda\cap (\lambda+H^{(\tau,0,1)})$. 
\end{rem}

\begin{ex}\label{u5}
It was shown in~\cite[Theorem 2.56 and Example 2.57]{H2} (see also~\cite[Theorem 16 and Example 3]{H}) that the convex subsets of cyclotomic model
sets with underlying
$\Z$-module $\Z[\zeta_5]$ are determined
by the $X$-rays in the $\Z[\zeta_5]$-directions of $U_5:=\{o/\vert o\vert\,|\,o\in O\}\subset\mathbb{S}^1$, where $O:=\{(1+\tau)+\zeta_{5}$, $(\tau-1)+\zeta_{5}$, $-\tau+\zeta_{5},2\tau-\zeta_{5}\}\subset\Z[\zeta_5]\setminus\{0\}$. Consequently, as was shown in the proof of Theorem~\ref{main1ico}(a), the convex subsets of generic icosahedral model sets
$\varLambda$ with underlying
$\Z$-module $L$ are determined
by the $X$-rays in the 
$L^{(\tau,0,1)}$-directions of $U_{\rm ico}:=\Phi^{-1}(U_5)\subset\mathbb{S}^2$. 
\end{ex}

\begin{rem}\label{rempl}
Since, by
 the work of Pleasants~\cite{PABP2}, the $\Z[\zeta_5]$-directions of $U_5$ are well suited in order to yield dense lines in cyclotomic model
sets with underlying
$\Z$-module $\Z[\zeta_5]$, it follows that the set of
$L^{(\tau,0,1)}$-directions $U_{\rm ico}$ from Example~\ref{u5} is well suited in order to yield dense
lines in the corresponding slices $\varLambda\cap (\lambda+H^{(\tau,0,1)})$, $\lambda\in \varLambda$, of generic icosahedral model
sets $\varLambda$ with underlying
$\Z$-module $L$. In fact, these directions even yield dense lines in icosahedral
model sets $\varLambda$ as a whole; {\it cf.}~\cite{PABP2}. In particular, neighbouring slices of the
form $\varLambda\cap (\lambda+H^{(\tau,0,1)})$, $\lambda\in \varLambda$, are densely
occupied and hence well
separated. Consequently, neighbouring lines in any of the directions of
$U$ that meet at least one point of a fixed icosahedral
model set $\varLambda$ are sufficiently separated. It follows that,
 in the practice of quantitative HRTEM,
 the resolution coming from the above directions
 is likely to be rather high, which makes Theorem~\ref{main1ico} look
 promising. 
\end{rem}

Finally, we want to demonstrate that, in an approximative sense, part (a) of Theorem~\ref{main1ico} even holds in the non-anchored case for regular generic icosahedral model sets. Before, we need a consequence of Weyl's theory of
uniform distribution; {\it cf.}~\cite{Weyl}. This analytical property of
regular icosahedral model sets was analyzed in general
in~\cite{Schl2}, \cite{Schl} and \cite{Moody2}. We need the following variant which relates the centroids of images of certain finite subsets of a regular icosahedral model set $\varLambda$ under the star map to the centroid of its window.

\begin{theorem}\label{weyl}
Let $\varLambda$ be a regular icosahedral model set of the form $\varLambda=\varLambda_{\rm ico}(0,W)$. Then, for
all $a\in\R^3$, one has the identity
$$
\lim_{R\rightarrow \infty}\,\, \frac{1}{\operatorname{card}(\varLambda\cap
  B_{R}(a))}\sum_{\alpha\in \varLambda\cap
  B_{R}(a)}\alpha^{\star}\,\,=\,\,\frac{1}{\operatorname{vol}(W)}\int_{W}y\,{\rm d}\lambda(y)\,,
$$
where $\lambda$ denotes the Lebesgue measure on $\R^3$.
\end{theorem}
\begin{proof}
This is a consequence of the uniform distribution
of the points of $\varLambda^*$ in the window, which gives
the integral by Weyl's lemma. The proof of the uniform
distribution property for model sets can be found in~\cite{Schl2,Moody,Moody2}.
\end{proof}

The following properties of sets $U\subset \mathbb{S}^1$ consisting of four pairwise non-parallel
$\Z[\zeta_5]$-directions will be of crucial importance:
\begin{itemize}
\item[(C)]
For all cyclotomic model sets
$\varLambda_{\rm cyc}$ with underlying
$\Z$-module $\Z[\zeta_5]$, the set
$\mathcal{C}(\varLambda_{\rm cyc})$ is determined by the $X$-rays in
the directions of $U$.
\item[(E)]
$U$ contains two directions of the form $o/ \vert o\vert,o'/ \vert o'\vert$, where $o,o'\in\Z[\zeta_5]\setminus\{0\}$ satisfy the relation $$\alpha_o\beta_{o'}-\beta_{o}\alpha_{o'}\,\,\in\,\,\Z[\tau]^{\times}=\{\tau^s\,|\,s\in\Z\}\,,$$ where the elements $\alpha_o,\alpha_{o'},\beta_{o},\beta_{o,}\in\Z[\tau]$ are determined by $o=\alpha_{o}+\beta_{o}\zeta_{5}$ and
$o'=\alpha_{o'}+\beta_{o'}\zeta_{5}$; {\it cf.} Remark~\ref{remzm2}. 
\end{itemize}

\begin{ex}\label{u5b}
The set $U_5\subset \mathbb{S}^1$ of four pairwise non-parallel
$\Z[\zeta_5]$-directions as defined in Example~\ref{u5} has property (C) by~\cite[Example 2.57]{H2} (see also~\cite[Example 3]{H}). Additionally, one can easily see that $U_5$ also has property (E).
\end{ex}

The significance of property (E) is expressed by the following result.

\begin{prop}\label{oneequiv}
Let $U\subset \mathbb{S}^1$ be a set of four pairwise
non-parallel $\Z[\zeta_5]$-directions with property {\rm (E)}. Then, for all finite subsets $F$ of $\Z[\zeta_5]$, one
has the inclusion $$G^{F}_{U}\,\,\subset\,\, \Z[\zeta_5]\,.$$
\end{prop}
\begin{proof}
This follows from~\cite[Theorem 1.130]{H2} (see also~\cite[Theorem 12]{H}).
\end{proof}

We are now able to show that, in an approximative sense to be
clarified below, for any fixed window $W\subset\R^3$ whose boundary $\operatorname{bd}(W)$ has Lebesgue
measure $0$ in $\R^3$, the set 
$\cup_{\varLambda\in
   \mathcal{I}^{\rm B}_{g}(W)}\mathcal{C}(\varLambda)$ (resp., $\cup_{\varLambda\in
   \mathcal{I}^{\rm F}_{g}(W)}\mathcal{C}(\varLambda)$)
 is determined by the $X$-rays in any set of
 $L^{(\tau,0,1)}$-directions, where $L=\operatorname{Im}(\Ic)$ (resp.,
 $L=\Ic_0$), of the form $U:=\Phi^{-1}(U')$, where $U'$ is a set of four pairwise non-parallel
$\Z[\zeta_5]$-directions with the properties {\rm (C)} and {\rm
   (E)}. Since the arguments for the F-type case and the B-type case are similar, we present the details for the B-type case only. Let $$F,F'\,\,\in\,\, \bigcup_{\varLambda\in
   \mathcal{I}^{\rm B}_{g}(W)}\mathcal{C}(\varLambda)\,,$$ say $F\in\mathcal{C}(\varLambda^{\rm B}_{\rm ico}(t,s+W))$ and $F'\in\mathcal{C}(\varLambda^{\rm B}_{\rm ico}(t',s'+W))$, where $t,t',s,s'\in\R^3$, and suppose that $F$ and $F'$ have the same $X$-rays in the directions of $U$. If $F=\varnothing$, then, by Lemma~\ref{cardinality}(a), one also gets
$F'=\varnothing$. One may thus assume, without loss of generality, that $F$ and $F'$ are
non-empty. Hence, there is an element $\lambda\in F$ such that $F\cap(\lambda+H^{(\tau,0,1)})$ and $F'\cap(\lambda+H^{(\tau,0,1)})$ are non-empty finite sets with the same $X$-rays in the directions of $U$. Then, by Lemma~\ref{trick7}, the non-empty finite subset $\Phi((F\cap(\lambda+H^{(\tau,0,1)}))-\lambda)$ of $\Z[\zeta_5]$ ({\it cf.} Lemma~\ref{res}) and the non-empty finite subset $\Phi((F'\cap(\lambda+H^{(\tau,0,1)}))-\lambda)$ of $\C$ have the same $X$-rays in the four pairwise non-parallel $\Z[\zeta_5]$-directions of $\Phi(U)=U'$. Then, by Lemma~\ref{fgrid} and Proposition~\ref{oneequiv} in conjunction with property~(E), one obtains
$$%\label{fton225}
\Phi((F\cap(\lambda+H^{(\tau,0,1)}))-\lambda), \Phi((F'\cap(\lambda+H^{(\tau,0,1)}))-\lambda)\subset G^{\Phi((F\cap(\lambda+H^{(\tau,0,1)}))-\lambda)}_{U'}\subset\Z[\zeta_5]\,.
$$
Thus, one gets
\begin{equation}\label{fton225}
F\cap(\lambda+H^{(\tau,0,1)}),F'\cap(\lambda+H^{(\tau,0,1)})\,\,\subset \,\,t+L\,.
\end{equation}
Since $F'\cap(\lambda+H^{(\tau,0,1)})\subset t'+L$, 
Relation~(\ref{fton225}) implies that $t+L$ meets $t'+L$, the latter
being equivalent to the identity $t+L=t'+L$. Note also that 
 the identity $t+L=t'+L$ is equivalent to the relation $t'-t\in L$. Clearly, one
has $$F-t\,\,\in\,\,\mathcal{C}(\varLambda^{\rm B}_{\rm ico}(0,s+W))\,.$$
Moreover, since the equality
$$\varLambda^{\rm B}_{\rm ico}(t'-t,s'+W)=\varLambda^{\rm B}_{\rm ico}(0,(s'+(t'-t)^{\star})+W)$$
holds, one further obtains 
$$F'-t\,\,\in\,\,\mathcal{C}(\varLambda^{\rm B}_{\rm ico}(t'-t,s'+W))\,\,=\,\,\mathcal{C}(\varLambda^{\rm B}_{\rm ico}(0,(s'+(t'-t)^{\star})+W))\,.$$
Clearly, $F-t$ and $F'-t$ again have the same $X$-rays in the
directions of $U$. Hence, by Lemma~\ref{cardinality}(b), $F-t$ and $F'-t$
have the same centroid. Since the star map $.^{\star}$ is $\Q$-linear, it follows that the finite subsets
$(F-t)^{\star}$ and $(F'-t)^{\star}$ of
$\R^3$ also have the same
centroid. Now, if one has $$F-t\,\,=\,\,B_{R}(a)\,\,\cap\,\,
\varLambda^{\rm B}_{\rm ico}(0,s+W)$$ and
$$F'-t\,\,=\,\,B_{R'}(a')\,\,\cap\,\,
\varLambda^{\rm B}_{\rm ico}(0,(s'+(t'-t)^{\star})+W)$$ for
suitable $a,a'\in\R^3$ and large $R,R'>0$ (which is rather natural in
practice), then Theorem~\ref{weyl} allows us to write

\begin{eqnarray*}
\frac{1}{\operatorname{vol}(W)}\int_{s +W}y\,{\rm
  d}\lambda(y)&\approx&
\frac{1}{\operatorname{card}\left(F-t\right)}\sum_{x\in
  F-t}x^{\star}\\ &=&
\frac{1}{\operatorname{card}\left(F'-t\right)}\sum_{x\in
  F'-t}x^{\star}\\
&\approx&\frac{1}{\operatorname{vol}(W)}\int_{(s'+(t'-t)^{\star})
  +W}y\,{\rm d}\lambda(y)\,.
\end{eqnarray*}
Consequently,  
$$
s+\int_{W}y\,{\rm
  d}\lambda(y)\,\,\approx\,\,(s'+(t'-t)^{\star})+\int_{W}y\,{\rm d}\lambda(y)\,,
$$
and hence $s\approx s'+(t'-t)^{\star}$. The latter means
that, approximately, both
$F-t$ and $F'-t$ are elements of the set
$\mathcal{C}(\varLambda^{\rm B}_{\rm ico}(0,s+W))$. Now, it follows
in this approximative sense 
 from property~(C) and Theorem~\ref{main1ico} that $F-t\approx F'-t$,
 and, finally, $F\approx F'$.

\begin{rem}
The above analysis suggests that, for all fixed windows $W\subset\R^3$ whose boundary $\operatorname{bd}(W)$ has Lebesgue measure $0$ in $\R^3$, the sets of the form $\cup_{\varLambda\in
   \mathcal{I}^{\rm B}_{g}(W)}\mathcal{C}(\varLambda)$ (resp., $\cup_{\varLambda\in
   \mathcal{I}^{\rm F}_{g}(W)}\mathcal{C}(\varLambda)$) are approximately determined by the $X$-rays in the four prescribed $L^{(\tau,0,1)}$-directions of $U_{\rm ico}$, where $L=\operatorname{Im}(\Ic)$ (resp.,
 $L=\Ic_0$); {\it cf.} Examples 6.12 and~6.15. Additionally, in the practice of quantitative HRTEM, the resolution coming from the directions of $U_{\rm ico}$ is likely to be rather high, which makes this approximative result look even more promising in view of real applications; {\it cf.} Remark~\ref{rempl}.
\end{rem}

\section{Outlook}
For a more extensive account of both uniqueness and computational complexity results in the discrete tomography of Delone sets with long-range order, we refer the reader to~\cite{H2}. This reference also contains results on the interactive concept of \emph{successive determination} of finite sets by $X$-rays and further extensions of settings and results that are beyond our scope here; compare also~\cite{H}. Although the results of this text and of~\cite{H2} give satisfying answers to the basic problems of
discrete tomography of icosahedral model sets, there is still a lot to
do to create a tool that is as satisfactory for the
application in materials science as is computerized tomography in its
medical or other applications. First, we believe
that it is an interesting problem to characterize the sets of
$\varLambda$-directions {\em in general position} having
the property that, for all icosahedral model sets
$\varLambda$, the set of convex subsets of $\varLambda$ is determined by the $X$-rays in
these directions; compare~\cite[Problems 2.1 and 2.3]{G}. Secondly, it
would be interesting to have experimental tests in order to see how well the above results work in
practice. Since there is always some noise involved when
physical measurements are taken, the latter also requires the ability to work with imprecise
data. For this, it is necessary to study stability and instability results in the discrete
tomography of icosahedral model sets in the future; {\it cf.}~\cite{AG}.

\section*{Acknowledgements}
It is my pleasure to thank Michael Baake, Uwe Grimm, Peter Gritzmann, Barbara Langfeld and Reinhard L\"uck for valuable discussions and
suggestions. This work was supported by the German Research Council (DFG), within the CRC 701, and by EPSRC, via Grant EP/D058465/1.


\begin{thebibliography}{12}
\bibitem{AG} Alpers A and Gritzmann P 2006 On stability, error correction, and noise compensation in discrete tomography {\em SIAM J. Discrete Math.} \textbf{20} 227--239

\bibitem{baake} Baake M 1997 Solution of the coincidence problem in
     dimensions $d\leq 4$. In: {\em The Mathematics
   of Long-Range Aperiodic Order} (Ed. Moody R V) pp. 9--44. (NATO-ASI Series C {\bf 489}, Kluwer,
 Dordrecht) Revised version \texttt{arXiv:math/0605222v1 [math.MG]
}

\bibitem{B} Baake M 2002 A guide to mathematical quasicrystals.
 In: {\em
    Quasicrystals. An Introduction to Structure, Physical Properties,
    and Applications} (Eds. Suck J-B, Schreiber M and H\"aussler
  P) pp. 17--48. (Springer, Berlin) \texttt{arXiv:math-ph/9901014v1}

\bibitem{BG2} Baake M, Gritzmann P, Huck C, Langfeld B and Lord
  K 2006 Discrete tomography of planar model sets 
  A{\bfseries 62} 419--433 \texttt{arxiv:math.MG/0609393} 

\bibitem{BH} Baake M and Huck C 2007 Discrete tomography of Penrose
    model sets {\em Philos. Mag.} {\bf 87} 2839--2846 \texttt{arXiv:math-ph/0610056v1}

\bibitem{BM} Baake M and Moody R V (Eds.) 2000 {\em Directions in Mathematical Quasicrystals} (CRM Monograph Series, vol. {\bf 13}, AMS, Providence, RI)

\bibitem{BPR} Baake M, Pleasants P A B and Rehmann U 2006 Coincidence
  site modules in $3$-space {\em Discr. Comput. Geom.} {\bf 38} 111--138 \texttt{arXiv:math/0609793v1 [math.MG]}

\bibitem{cmp} Chen L, Moody R V and Patera J (1998)
  Non-crystallographic root systems. In: {\em Quasicrystals and
    Discrete Geometry} (Ed. Patera J) pp. 135--178. (Fields Institute
  Monographs, vol. 10, AMS, Providence, RI) 

\bibitem{con} Conway J H and Sloane N J A 1999 {\em Sphere Packings,
Lattices and Groups} (3rd ed., Springer, New York) 

\bibitem{cow} Cowley J M 1995 {\em Diffraction Physics}
 (3rd rev. ed., North-Holland, Amsterdam)


\bibitem{d} Dubois J-M 2005 {\em Useful Quasicrystals} (World Scientific, Singapore)

\bibitem{few} Fewster P F 2003 {\em X-ray Scattering from
    Semiconductors} (2nd ed., Imperial College Press, London)

\bibitem{G} Gardner R J 2006 {\em Geometric Tomography} (2nd ed., Cambridge University Press, New York)

\bibitem{GG} Gardner R J and Gritzmann P 1997 Discrete tomography:
  determination of finite sets by X-rays {\em Trans. Amer. Math. Soc.} {\bf 349} 2271--2295

\bibitem{GG2} Gardner R J and Gritzmann P 1999 Uniqueness and
    complexity in discrete tomography. In:~\cite{HK} pp. 85--114

\bibitem{GGP} 
  Gardner R J, Gritzmann P and Prangenberg D 1999  
  On the computational complexity of
  reconstructing lattice sets from their X-rays 
  {\em Discrete Math.} \textbf{202} 45--71

\bibitem{Gr}
Gritzmann P 1997 On the reconstruction of finite lattice sets
from their X-rays. In: {\em Lecture Notes on Computer Science} (Eds.:
Ahronovitz E and Fiorio C) pp. 19--32 (Springer, London)

\bibitem{GL} Gritzmann P and Langfeld B On the index of Siegel grids and its application to the tomography of quasicrystals \emph{Eur. J. Comb.} to appear

\bibitem{gu} Guinier A 1994 {\em X-ray Diffraction in Crystals,
    Imperfect Crystals, and Amorphous Bodies} (Dover Publications, New
  York)

%\bibitem{gkb} Guyot P, Kramer P and de Boissieu M 1989 Quasicrystals
%  \emph{Rep. Prog. Phys.} {\bf 79} 1373--1425 

\bibitem{hw} Hardy G H and Wright E M 1979 {\em An Introduction to the
  Theory of Numbers} (5th ed., Clarendon Press, Oxford)

\bibitem{HK} Herman G T and Kuba A (Eds.) 1999 
{\em Discrete Tomography: Foundations, Algorithms, and Applications}
(Birkh\"auser, Boston)

\bibitem{HG} de Boissieu M, Guyot P and Audier M 1994 Quasicrystals:
  quasicrystalline order, atomic structure and phase transitions. In:
  {\em Lectures on Quasicrystals} (Eds. Hippert F and Gratias D) pp. 1--152 (EDP Sciences, Les Ulis)

\bibitem{H} Huck C 2007 Uniqueness in discrete tomography
  of planar model sets. Notes \texttt{arXiv:math/0701141v2 [math.MG]}

%\bibitem{H1} Huck C 2007 Uniqueness in discrete tomography of Delone sets with long-range order. Submitted \texttt{arXiv:0711.4525v1 [math.MG]}

\bibitem{H2} Huck C 2007 {\em Discrete Tomography of Delone Sets with Long-Range Order} (Ph.D. thesis (Universit\"at Bielefeld), Logos Verlag, Berlin)


\bibitem{j} Janot C 1994
{\em Quasicrystals: A Primer} (2nd ed., Clarendon Press, Oxford)

\bibitem{ks} Kisielowski C, Schwander P, Baumann F H, Seibt M, Kim Y
  and Ourmazd A 1995 An approach to quantitative high-resolution
  transmission electron microscopy of crystalline materials {\em Ultramicroscopy} \textbf{58} 131--155

\bibitem{Moody} Moody R V 2000 Model sets: a survey. In: {\em From
    Quasicrystals  to More Complex Systems} (Eds. Axel F,
  D\'{e}noyer F and Gazeau J-P) pp. 145--166 (EDP Sciences, Les Ulis,
  and Springer, Berlin) \texttt{arXiv:math/0002020v1 [math.MG]} 

\bibitem{Moody2} Moody R V 2002 Uniform distribution in model
    sets {\em Canad. Math. Bull.} \textbf{45} 123--130

\bibitem{MP} Moody R V and Patera J 1993 Quasicrystals and icosians  {\em J. Phys. A} \textbf{26} 2829--2853 

\bibitem{PABP} Pleasants P A B  2000 Designer quasicrystals:
  cut-and-project sets with pre-assigned properties. In: ~\cite{BM}
  pp. 95--141

\bibitem{PABP2} Pleasants P A B 2003 Lines and planes in $2$- and
  $3$-dimensional quasicrystals. In: {\em Coverings of Discrete
    Quasiperiodic Sets} (Eds. Kramer P and Papadopolos Z) pp. 185--225 (Springer Tracts in Modern Physics, vol. \textbf{180}, Springer, Berlin)

%\bibitem{ps} Preparata F P and Shamos M I 1985 {\em Computational Geometry: An Introduction}
%  (Springer, New York)

\bibitem{Re} R\'enyi A 1952 On projections of probability
    distributions {\em Acta Math. Acad. Sci. Hung.} \textbf{3} 131--142


\bibitem{s} Shechtman D, Blech I, Gratias D and Cahn J W 1984
  Metallic phase with long-range orientational order and no
  translational symmetry {\em Phys. Rev. Lett.} \textbf{53} 1951--1953

\bibitem{Schl2} Schlottmann M 1998 Cut-and-project sets in locally
  compact Abelian groups. In: {\em Quasicrystals and Discrete
    Geometry} pp. 247--264 (Fields Institute Monographs, vol. {\bf 10}, AMS, Providence, RI)

\bibitem{Schl} Schlottmann M  2000 Generalized model sets and dynamical
  systems. In: ~\cite{BM} pp. 143--159 

\bibitem{sk} Schwander P, Kisielowski C, Seibt M, Baumann F H, Kim Y
  and Ourmazd A 1993 Mapping projected potential, interfacial
  roughness, and composition in general crystalline solids by
  quantitative transmission electron microscopy {\em Phys. Rev. Lett.} \textbf{71} 4150--4153


\bibitem{so} Steinhardt P J and Ostlund S (Eds.) 1987 
{\em The Physics of Quasicrystals} 
 (World Scientific, Singapore)

%\bibitem{schw} Schwarzenberger R L E 1980 {\em N-dimensional
%    crystallography} (Pitman, San Francisco) 


\bibitem{St} Steurer W 2004 Twenty years of structure research on
  quasicrystals. Part I. Pentagonal, octagonal, decagonal and
  dodecagonal quasicrystals {\em Z. Kristallogr.} {\bf 219} 391--446 

\bibitem{Wa} Washington L C 1997 {\em Introduction to Cyclotomic Fields} (2nd ed., Springer, New York)

\bibitem{Weyl} Weyl H 1916 \"Uber die Gleichverteilung von Zahlen
    mod. Eins {\em Math. Ann.} \textbf{77} 313--352

\end{thebibliography}
\end{document}